\documentclass{elsarticle}
\usepackage{lineno,hyperref}
\usepackage{ulem}
\usepackage{lineno,hyperref}
\usepackage{amsmath,amssymb}
\usepackage{graphicx}
\usepackage{amsmath}
\usepackage{multirow}
\usepackage{color}
\usepackage{makecell}
\usepackage{amsmath,amssymb}
\usepackage{diagbox}
\biboptions{numbers,sort&compress}
\usepackage{subfigure}% Support for small, `sub' figures and tables
\usepackage{slashbox}
\usepackage{float}

\newtheorem{thm}{Theorem}
\newtheorem{lem}[thm]{Lemma}
\newdefinition{rmk}{Remark}
\newproof{pf}{Proof}

\allowdisplaybreaks

\journal{Journal of \LaTeX\ Templates}

\usepackage{numcompress}

\begin{document}

\begin{frontmatter}

\title{On a progressive and iterative approximation method with memory for least square fitting \tnoteref{mytitlenote}}
\tnotetext[mytitlenote]{Fully documented templates are available in the elsarticle package on \href{http://www.ctan.org/tex-archive/macros/latex/contrib/elsarticle}{CTAN}.}

\author[1]{Zheng-Da Huang\corref{mycorrespondingauthor}}
\cortext[mycorrespondingauthor]{Corresponding author}
\ead{zdhuang@zju.edu.cn}

\author[1]{Hui-Di Wang}\ead{hdwang@cjlu.edu.cn}%\fnref{myfootnote}}
%\address{Hangzhou, P.R. China}

\address[1]{School of Mathematical Sciences, Zhejiang University, Hangzhou 310027, P.R. China}
\address[2]{College of Sciences, China Jiliang University, Hangzhou 310018, P.R. China}

%\author[mymainaddress,mysecondaryaddress]{Hong-Wei Lin}
%\author[mysecondaryaddress]{Zheng-Da Huang\corref{mycorrespondingauthor}}

\begin{abstract}
In this paper, we present a progressive and iterative approximation method with memory for least square fitting
 (MLSPIA). It adjusts the control points  and the weighted sums iteratively to construct a series
of fitting curves (surfaces) with three weights. For any normalized totally positive basis even when the collocation matrix is of  deficient column rank, we obtain a condition to guarantee that these curves (surfaces) converge to the least square fitting curve (surface) to the given data points. It is proved that the theoretical convergence rate of the method is faster than the one of the progressive and iterative approximation method for least square fitting (LSPIA) in [Deng C-Y, Lin H-W. Progressive and iterative approximation for least squares B-spline curve and surface fitting. Computer-Aided Design 2014;47:32-44] under the same assumption. Examples verify this phenomenon.
\end{abstract}


\end{frontmatter}

%\linenumbers

\section{Introduction}

Least square fitting to a set of data points by a parametric curve or surface is a fundamental problem in Computer Aided Geometric Design \cite{RN2222QI,Farin} and a well-known key technology in many application areas. It is powerful  especially when dada set is very larger.

Based on different metric descriptions, there are two main kinds of least square fittings. One is to take control points and parameter knots as variables, so that nonlinear least square problems  should be solved. Since it is nonlinear, a proper initial fitting need. The other one is to assign a parameter for each data point, so that the data set becomes organized. In this case, linear least square fittings need to be solved to obtain control points. Since it is linear, the construction of the method may be simpler. One can look for more details about these two kinds of fittings, for examples, in \cite{Borgescf,RN402,RN214GI,RN209Tspline,RN2222QI,RN183LSPIA,GalvezaA,RN181TwoWeights,RN209singular,Ebrahimi,RN230Survey,Farin,
RN220LIU,Vaitkus} and references therein.

 To improve the smoothness of fitting curves or surfaces, fairness terms with user-defined parameters may be considered sometimes\cite{Farin,Vaitkus,Ebrahimi,RN402,RN214GI}. As you know, choices of these user-defined parameters may affect the smoothness and accuracy, but there is no optimal choice of these user-defined parameters. %There are many literatures published around both of these two kinds of fittings.

To obtain  better fitting efficiency, researchers also consider the way of how to determine a relatively low number of control points used in the least square fitting \cite{GalvezaA}.% and how to assign parameters to points in data sets\cite{Vaitkus}.

In this paper, we consider the linear least square fitting, with a given number of control points, to the data set of which each point is assigned to a parameter value. For simplicity, we also assume that there are no fairness terms.

Standard fittings to any data set of such kind are to obtain control points by solving corresponding systems of linear equations directly. When occasionally extra control points need to be added to improve the fitting accuracy, new systems of linear equations should be assembled and resolved. The LSPIA method, a least square version of PIA method\cite{RN218Totally}, developed in \cite{RN183LSPIA}, is one of the least square fitting methods that need not solve any system of linear equations directly.  It approximates control points iteratively, and generates a series of curves (surfaces) converging to the least square fitting result of the given data points. It is easy to make the fitting result hold the shape preserving
property. When an extra control point should be added, rather than assembling a new system of linear equations, only iterative formula needs to be modified, as said in \cite{RN183LSPIA,RN230Survey}.

To increase the convergence speed, in this paper, a PIA method with memory  for least square fitting (abbr. MLSPIA) is presented.  As usual, the method constructs a series of curves (surfaces) by adjusting the control points and the weighted sums iteratively. Reasons for us to call the  new method as MLSPIA method can be found in Section 3.

Our contributions to the MLSPIA method are:
\begin{itemize}
\item We prove in a same way that the MLSPIA method is feasible regardless of whether the collocation matrix is of deficient column rank or not.
\item  Compared with the LSPIA method in \cite{RN183LSPIA}, the MLSPIA method does not increase too  high cost, especially for the case of a  very large data set, while the convergence rate is much faster when the collocation matrix is of  full column rank. %in \cite{RN209Tspline}
\item The MLSPIA method preserves the advantages of the LSPIA method, including handling point set of large size and allowing the adjustment of the number
of control points and knot vectors in the iterations, the shape preserving property and parallel performance.
\end{itemize}

The paper is organized as follows. In Section 2, related works are introduced briefly, and in Section 3, the iterative format and the convergence analysis of the MLSPIA method are presented. We compare the convergence speed of the MLSPIA method with that of the LSPIA method in Section 4, and propose the MLSPIA iterative format of the surface in Section 5. Finally, five representative  examples are performed in Section 6.

\section{Related work}

PIA methods, i.e., geometric fitting methods that make use of the progressive and iterative approximation property of univariate NTP bases, have the similar iterative formats as the geometric interpolation (GI) methods, such as in \cite{RN191,RN212GI}, etc. The PIA methods depend on the parametric distance while the GI methods rely on the geometric distance. The study starts from the woks of de Boor \cite{RN209de}, Qi and co-authors \cite{RN208QI}, and Yamaguji \cite{RN1977Y} in 1970's. Rather than those  methods based directly on the solution of  a linear system, by fixing parameters, each PIA method generates a series of curves (or surfaces) to approximate the interpolation curve (or surface).

Give a point set $\{Q_{i}\}_{i=1}^m$ and %select $\{P_{i}\}_{i=1}^n$ from $\{Q_{i}\}_{i=1}^m$ with
$m \ge n$. Let $0\le t_1<t_2< \cdots <t_m \le 1$ be an increasing sequence.

Choose $ B_1 (t), B_2(t), \cdots, B_n(t)$ as an normalized totally positive (NTP) basis of a vector space of real functions defined on $[0,1]$. As you know, the basis is said to be NTP if $B_i(t) \ge 0, i\!=\!1, 2, \cdots, n,$ and  $\sum\limits_{i=1}^n B_i(t)=1$  hold for all $t\in [0, 1]$, and its collocation matrix at any increasing sequence is a totally positive matrix. A matrix is called to be totally positive if all of its minors are nonnegative \cite{RN225Convexity,RN226Delgado}.

Clearly, the collocation matrix of the basis $B_1 (t), B_2(t), \cdots, B_n(t)$ at the given increasing sequence $0\le t_1<t_2< \cdots <t_m \le 1$ is
\[B_{m \times n}=
\begin{pmatrix}
B_1(t_1)\ &B_2(t_1)\ &\cdots \ &B_n(t_1)\\
B_1(t_2)\ &B_2(t_2)\ &\cdots \ &B_n(t_2)\\
\vdots &\vdots &\ddots &\vdots \\
B_1(t_m)\ &B_2(t_m)\ &\cdots \ &B_n(t_m)
\end{pmatrix}.\]

 Choosing  any  $\left\{P_i^0\right\}^n_{i=1}$ in the vector space, where the set $\left\{Q_i\right\}^m_{i=1}$ is located,  as the initial control points set, the PIA method defined in \cite{RN218Totally} (with $m=n$) and the LSPIA method in \cite{RN183LSPIA} (with $m>n$) approximate the target curve with the NTP basis $B_1(t), B_2(t), \cdots, B_n(t)$ by the following curves iteratively:
\begin{equation}\label{curve}
C^{k+1}(t)=\sum\limits_{i=1}^n B_i (t) P_i^{k+1},\ \ t \in [0,1],\ \ \ \ k \ge -1,
\end{equation}
with
\begin{equation}\label{PIA}
 P_i^{k+1}=P_i^k+ \Delta_i^k,\ \ \ \ i=1,2,\ \cdots, n,\ \ \ \ k \ge 0,
\end{equation}
where %for each $k \ge 0$,
\[
\Delta_i^k \!=\!
\left\{\!\! \begin{array}{cr}
\!Q_i\!-\!C^k(t_i), \!&\mbox{\small \!for PIA method \cite{RN218Totally}},\\
\!\mu \sum\limits_{j\!=\!1}^m B_i(\!t_j\!)\!\left(\!Q_j\!-\!C^k(t_j)\!\right)\!,\! &\mbox{\small \!for LSPIA method \cite{RN183LSPIA}},
\end{array}
\right.
\ i=1, 2, \cdots\!,\! n, \ \ k \ge 0,
\]
and $\mu > 0$ is a constant.

The PIA method in the format of \eqref{curve} and \eqref{PIA}, defined and firstly named in \cite{RN199Constructing}, is based on the idea of profit and loss modification proposed in \cite{RN208QI,RN209de} by using the non-uniform B-spline basis. It is generalized to the method in the same format as in \cite{RN218Totally} by replacing the non-uniform B-spline basis with any totally positive basis whose collocation matrix is nonsingular.

The PIA method in \cite{RN218Totally} has been developed to a weighted PIA method in \cite{RN193Weighted} that  speeds up the convergence, and  to an extended PIA (EPIA) format with NTP bases in \cite{RN188extended} that focuses on large-scale data sets. It is also extended to approximate the NURBS curves (surfaces) in \cite{RN211nurbs} by using homogeneous coordinates.

In \cite{RN192Local}, the author considers the local property of the progressive iteration approximation, and proposes a PIA method with local format. It is shown that the limit curve will still interpolate the subset of  data points if only a subset of the control points are adjusted iteratively, and the others remain unchanged.

The PIA method in \cite{RN218Totally} has been  modified to approximate surface by constructing  series of tensor product surfaces based on the corresponding tensor bases to fitting in \cite{RN199Constructing,RN218Totally,RN193Weighted,RN188extended,RN192Local},
to approach the triangular B\'ezier surface or the rational triangular B\'ezier surface in \cite{RN189triangular,RN221HU} by using the  non-tensor product triangular Bernstein basis. It is also generalized to construct the progressive interpolation (PI) methods for subdivision surface fittings such as the Loop subdivision fitting \cite{RN202Loop,RN201WeightedPI}, the Catmull-Clark subdivision fitting \cite{RN207Catmull-Clark} and the Doo-Sabin subdivision surface fitting \cite{RN210shape}.

The PIA methods converge with suitable weights, and  own the advantages of convexity preserving and the explicit expression of curves (or surfaces), locality and adaptivity. The locality and adaptivity makes it possible  for us to improve the approximation precise to a segment of the result curve (or surface) by changing only control points that link with this segment with lower cost \cite{RN199Constructing,RN218Totally}. One can look for more details in the survey paper \cite{RN230Survey}.

The LSPIA method in the format of \eqref{curve} and \eqref{PIA}, defined in \cite{RN183LSPIA}, is constructed mainly for very large-scale  data sets,  and inherits advantages of PIA methods.  It is developed by using the T-splines in \cite{RN209TSPLINE}. By replacing the B-spline basis with the generalized B-spline basis, it is generalized in \cite{RN181TwoWeights} to the weighted least square fitting  curve, and, by replacing the tensor product B-spline basis with the non-tensor product bivariate B-spline basis,  extended in \cite{RN220LIU} to the regularized least square fitting surface.

The convergence of the LSPIA methods in the singular case is first proved in \cite{RN209singular} for the LSPIA method published in \cite{RN183LSPIA} and some variants of it. In \cite{RN183LSPIA}, only convergence in the  nonsingular cases is proved. The  way used in \cite{RN209singular} to prove the convergence in the singular case is different from that used in \cite{RN183LSPIA} to prove the convergence in the nonsingular case.

The history of a least square fitting to a curve with B-spline function can be traced back to  \cite{RN400} at least, where authors computed the least square cubic spline  fitting with fixed knots based on the Gram-Schmidt process.

\section{The MLSPIA method and its convergence}\label{MLSPIA_curves_sec}

In this section, we will propose the MLSPIA method, i.e., the progressive iteration approximation method with memory for least square fitting, and discuss  its convergence.

For the given  data point set $\{Q_i\}_{i=1}^m$ and  any initial control point set  $\{P_i^0\}_{i=1}^n $,  with $m>n$, in the vector space where  $\{Q_i\}_{i=1}^m$ is located, the NTP basis $B_1(t)$, $B_2(t)$, $\cdots, B_n(t)$ and the increasing sequence $0=t_1<t_2< \cdots < t_m=1$ , the MLSPIA method approximates the curve iteratively by \eqref{curve} and \eqref{PIA} with
\begin{equation}\label{Iter}
\left\{\!\!\!
\begin{array}{lc}
\Delta_i^{0}= \upsilon \sum \limits_{j = 1}^m B_i(t_j) \Lambda_j^0, &\\
\Delta_i^{k}=(1-\omega) \Delta_i^{k-1} + \gamma \delta_i^k + (\omega  - \gamma ) \delta_i^{k-1}, &k \ge 1,\\
\delta_i^k=\upsilon \sum \limits_{j = 1}^m B_i(t_j)(Q_j - C^k(t_j)),
&k \ge 0,
\end{array}\ i = 1, 2, \cdots, n,
\right.\!\!
\end{equation}
where $\Lambda_j^0$ is an arbitrary point in the vector space where $\{Q_i\}_{i=1}^m$  is located for each $j =1,2,\cdots, m$, and $\omega, \gamma, \upsilon$ are three real weights.

There are at least two reasons for us to call the method,  defined by \eqref{curve}, \eqref{PIA} and \eqref{Iter},  the MLSPIA method.  It is called as a LSPIA method since the construction of it is very  similar to the LSPIA  method in the format of \eqref{curve} and \eqref{PIA}, defined in \cite{RN183LSPIA}, expect the difference of $\Delta_i^k$ for each possible $i$ and $k$.  As we can see later, it will also converge to the least square fitting curve of the given data point set $\{Q_i\}_{i=1}^m$ for suitable weights $\omega, \upsilon, \gamma$.  And  it is called  a method with memory since when we compute $\Delta_i^k$ via \eqref{Iter}, we need to store and use $\Delta_i^{k-1}$ and $\delta_i^{k-1}$ that have been computed, and compute $\delta_i^k$ for each $k \ge 1$ and $i=1,2, \cdots, n$. Since the price of DRAM becomes more and more cheaper, storing data of the previous step becomes more and more cheaper, too.  It is a non-stationary method since it needs to store the information of previous steps.

Compared with the LSPIA method  in \cite{RN183LSPIA}, the MLSPIA method just needs two more additions and three more scalar multiplications for each $k \ge 1$ and each $1 \le i \le n$. So it does not increase too much cost for the case of a very large data set. As we can see in the next section, the MLSPIA method  has the faster convergence rate.

In the following, we at first give an  equivalent expression of the MLSPIA method defined by \eqref{curve}, \eqref{PIA} and \eqref{Iter}.
\begin{lem}
For any  initial control point set $\{P_i^0\}_{i=1}^n$, and any point set $\{\Lambda_i^0\}_{i=1}^m$, all the $P_i^k$, $i=1,2,\cdots,n$, $k \ge 0$, generated by \eqref{PIA} and \eqref{Iter}, can be rewritten as
\begin{equation}\label{definition_of_P_k+1}
P_i^{k+1}=P_i^{k}+\upsilon \sum\limits_{j=1}^m B_i(t_j) \Lambda_j^k ,\ \ i=1, 2, \cdots, n, \ \ k \ge 0,
\end{equation}
where
\begin{equation}\label{Lam}
\left\{
\begin{array}{c}
\Lambda_i^{k+1}=(1-\omega)\Lambda_i^k + \omega Q_i-\sum\limits_{j_1=1}^n B_{j_1}(t_i)\left[\omega P_{j_1}^k+\gamma \upsilon  \sum\limits_{j_2=1}^m B_{j_1}(t_{j_2})\Lambda_{j_2}^k  \right],\\
\ \ \ \ \ \ \ \ \ \ \ \ \ i=1, 2, \cdots, m, \ \ k \ge 1.
\end{array}\right.
\end{equation}
\end{lem}
\begin{pf}
 When $k=0$, it can be easily checked that \eqref{definition_of_P_k+1} holds by \eqref{PIA} and \eqref{Iter}.

Suppose that \eqref{definition_of_P_k+1} holds for some $k\ge0$. Since, for each $1 \le i \le n$,
\begin{eqnarray*}
%\begin{aligned}
\omega \delta_i^k +\gamma(\delta_i^{k+1}-\delta_i^k)\!\!\!\!&=&\!\!\!\! \upsilon\omega\sum \limits_{j = 1}^m\! B_i(t_j)(Q_j-C^k(t_j)) \\
 \!\!\!\!&&\!\!\!\!\qquad \qquad- \upsilon\gamma\! \sum \limits_{j = 1}^m\! B_i(t_j)\left( C^{k+1}(t_j)-C^k(t_j)\right)\\
 &\!\!\!\!=&\!\!\!\!\upsilon\omega\sum \limits_{j = 1}^m\! B_i(t_j)\left(Q_j-\sum \limits_{j_1=1}^n\! B_{j_1}(t_j)P_{j_1}^k\right)\\
 \!\!\!\!&&\!\!\!\!\qquad\qquad-\upsilon\gamma \!\sum \limits_{j = 1}^m \!B_i(t_j)\!\sum \limits_{j_1=1}^n\!B_{j_1}(t_j)\left( P_{j_1}^{k+1}-P_{j_1}^k\right)
%\end{aligned}
\end{eqnarray*}
holds by  \eqref{Iter}, we have by \eqref{PIA} and \eqref{Iter} for each $1 \le i \le n$ that
\begin{eqnarray*}
\Delta_i^{k+1}\!\!\!\!&=&\!\!\!\! (1-\omega)\Delta_i^k+ \omega \delta_i^k + \gamma (\delta_i^{k+1}-\delta_i^k)\\
\!\!\!\!&=&\!\!\!\! (1-\omega)(P_i^{k+1}-P_i^k) +\upsilon\omega\sum \limits_{j = 1}^m\! B_i(t_j)\left(Q_j-\sum \limits_{j_1=1}^n\! B_{j_1}(t_j)P_{j_1}^k\right)\\
 &&\qquad\qquad\qquad\qquad\ -\upsilon\gamma \!\sum \limits_{j = 1}^m \!B_i(t_j)\!\sum \limits_{j_1=1}^n\!B_{j_1}(t_j)\left( P_{j_1}^{k+1}-P_{j_1}^k\right)\\
\!\!\!\!&=&\!\!\!\! (1-\omega)\upsilon\! \sum\limits_{j=1}^m \! B_i(t_j)\Lambda_j^{k}+ \upsilon\omega\sum \limits_{j = 1}^m\! B_i(t_j)\left(Q_j-\sum \limits_{j_1=1}^n\! B_{j_1}(t_j)P_{j_1}^k\right)\\
 &&\qquad\qquad\qquad\qquad\ \ - \upsilon^2\gamma\!\sum \limits_{j = 1}^m \!B_i(t_j)\!\sum \limits_{j_1=1}^n\!B_{j_1}(t_j)\! \sum\limits_{{j_2}=1}^m \! B_{j_1}(t_{j_2})\Lambda_{j_2}^{k}\\
\!\!\!\!&=&\!\!\!\!\upsilon \sum\limits_{j=1}^m \! B_i(t_j)\Lambda_j^{k+1},
\end{eqnarray*}
if we let
\begin{equation*}
\Lambda_i^{k+1}\!= (1-\omega)\Lambda_i^k + \omega Q_i -\!\!\sum\limits_{j_1=1}^n \!B_{j_1}(t_i)\!\!\left[\omega P_{j_1}^k+\gamma \upsilon \! \!\sum\limits_{j_2=1}^m \!B_{j_1}(t_{j_2})\Lambda_{j_2}^k  \right]\!\!, i=1,2,\cdots,m.
\end{equation*}
It follows that
\begin{equation*}
P_i^{k+2}=P_i^{k+1}+\Delta_i^{k+1}=P_i^{k+1}+\upsilon \sum\limits_{j=1}^m \! B_i(t_j)\Lambda_j^{k+1}, i=1,2,\cdots,n,
\end{equation*}
which shows that \eqref{definition_of_P_k+1} holds for $k+1$.\\
\indent By the induction method, \eqref{definition_of_P_k+1} holds for any $k \ge 0$. The proof is completed.
\end{pf}

Denote
\begin{equation}\label{Vec}
\left\{\begin{array}{l}
\Lambda^k=[\Lambda_1^k, \Lambda_2^k, \cdots, \Lambda_m^k]^T,\\
P^k=[P_1^k, P_2^k, \cdots, P_n^k]^T,
\end{array}
\right.\ \  k \ge 0,
\end{equation}
and
\begin{equation}\label{B}
Q=[Q_1, Q_2, \cdots, Q_m]^T,\ \ B=\left(B_i(t_j)\right)_{m \times n},
\end{equation}
where $P_i^k$, $i=1,2,\cdots,n$, and $\Lambda_i^k$, $i=1,2,\cdots,m$, are defined by \eqref{definition_of_P_k+1} and \eqref{Lam}, respectively, for each $k\ge0$.
Then \eqref{definition_of_P_k+1} and \eqref{Lam} can be written as
\begin{equation}\label{DIEDAI2}
\left\{\begin{array}{l}
\Lambda^{k+1}=\left[(1-\omega)I_m-\gamma \upsilon BB^T\right]\Lambda^k+\omega(Q-BP^k),\\
P^{k+1}=P^k+\upsilon B^T\Lambda^k,
\end{array} \right.\ k \ge 0,
\end{equation}
or,
\begin{equation}\label{juzhen_xingshi}
\left(
\begin{array}{l}
\Lambda^{k+1} \\
P^{k+1}
 \end{array} \right) = H_{\omega, \gamma, \upsilon}
 \left(
 \begin{array}{l}
\Lambda^k\\
P^k
 \end{array} \right)+C_{\omega, \gamma, \upsilon},
\  \ \ k \ge 0,
\end{equation}
where
\begin{equation}\label{H}
H_{\omega, \gamma, \upsilon}=\left(
\begin{array}{*{20}{c}}
(1-\omega)I_m-\gamma \upsilon BB^T&-\omega B\\
\upsilon B^T&I_n
\end{array}
\right),\
C_{\omega, \gamma, \upsilon}=
\left(
\begin{array}{c}
\omega Q\\
0
\end{array}
\right).
\end{equation}

\eqref{juzhen_xingshi} can be regarded as the matrix iterative format of the MLSPIA method. So to consider the convergence of the curves generated by the MLSPIA method, we only need to discuss the convergence of $\displaystyle \binom{\Lambda^k}{P^k}$ generated by \eqref{juzhen_xingshi} for the given initial point $\displaystyle \binom{\Lambda^0}{P^0}$.
\begin{thm}\label{Theorem}
Let $m>n$ and $B$ be defined by \eqref{B}. If $rank(B)=r$ and $\sigma_1 \ge \sigma_2 \ge \cdots \ge \sigma_r >0$ are the singular values of the matrix $B$, then, for any initial control point set $\{P_i^0\}_{i=1}^n$ and any point set $\{\Lambda_i^0\}_{i=1}^m$, the series of the curves, generated by the MLSPIA method, defined by \eqref{curve}, \eqref{PIA} and \eqref{Iter}, with weights $\omega, \gamma, \upsilon$,  converge to the least square fitting curve of the given $\{Q_i\}_{i=1}^m$ if
\begin{equation}\label{Convergence_1}
0<\omega<2,\ \ \ \omega-\frac{\omega}{\sigma_1^2 \upsilon}<\gamma<\frac{\omega}{2}-\frac{\omega-2}{\sigma_1^2 \upsilon}, \ \ \  \upsilon>0.
\end{equation}
\end{thm}

To prove Theorem \ref{Theorem}, we at first rewrite the equivalent format of the iteration matrix $H_{\omega, \gamma, \upsilon}$ defined by \eqref{H}.

Suppose that the singular value decomposition of $B$ is
\begin{equation}\label{qiyi}
B=U\left(\begin{array} {*{20}{c}} \Sigma_r & 0 \\ 0 & 0 \end{array} \right) V^T,
\end{equation}
where $U \in R^{m \times m}$ and $V \in R^{n \times n}$ are orthogonal matrices and $\Sigma_r =diag(\sigma_1, $$\sigma_2,$ $ \cdots, \sigma_r)$ with
$\sigma_1 \ge \sigma _2 \ge \cdots \ge \sigma_r >0.$ Let
\begin{equation}\label{N}
W=\left( \begin{array}{*{20}{c}}
U & 0\\ 0 & V
\end{array}\right),
\end{equation}
then $W \in R^{(m+n) \times (m+n)}$ is an orthogonal matrix. By the equalities
\[ U^TBV=\left( \begin{array}{*{20}{c}}  \Sigma_r & 0 \\ 0 & 0 \end{array} \right) \in R^{m \times n},\ \ V^TB^TU=\left( \begin{array}{*{20}{c}}  \Sigma_r & 0 \\ 0 & 0 \end{array} \right) \in R^{n \times m}, \]  and \[ U^TBB^TU=\left( \begin{array}{*{20}{c}}  \Sigma_r^2 & 0 \\ 0 & 0 \end{array} \right) \in R^{m \times m}, \] it holds that
\begin{equation}\label{xiangsi}
W^T H_{\omega, \gamma, \upsilon}W = \left( \begin{array}{*{20}{c}} \hat H_{\omega, \gamma, \upsilon} & 0 \\ 0 & I_{n-r}  \end{array}\right),
\end{equation}
where
\begin{equation}\label{hat}
\hat H_{\omega, \gamma, \upsilon}=\left( \begin{array}{*{30}{c}}
 (1-\omega)I_r-\gamma \upsilon \Sigma_r^2 & 0 & -\omega \Sigma_r\\
 0 & (1-\omega)I_{m-r} & 0\\
 \upsilon \Sigma_r & 0 & I_r
 \end{array}\right).
\end{equation}

The following lemmas are useful in analyzing the convergence of the method.
\vskip-0.2cm\begin{lem}\label{1}
Under the assumption of Theorem \ref{Theorem}, $\lambda$ is an eigenvalue of $\hat H_{\omega, \gamma, \upsilon}$ defined by \eqref{hat} if and only if $\lambda = 1- \omega$ or $\lambda$ is a root of the one of the following equations
\begin{equation}\label{polynomial}
\lambda^2\!+\![\gamma \upsilon \sigma_i^2 \!-\!(2\!-\! \omega)]\lambda\! + \! \sigma_i^2\upsilon(\omega-\gamma)+1-\omega=0, \ \ \ i=1, 2, \cdots, r.
\end{equation}
\end{lem}
\begin{pf}
Since $m>n$, it follows from \eqref{hat} that
\[
\begin{array}{ll}
|\hat H_{\omega, \gamma, \upsilon}\! -\! \lambda I_{m+r}|\!\!\!&\!\!=\!\left|
\begin{array}{*{30}{c}}
(1\!-\!\omega\!-\!\lambda)I_r-\gamma \upsilon \Sigma_r^2 & 0 & \!-\!\omega \Sigma_r\\
0 & (1\!-\!\omega\! -\!\lambda)I_{m-r} & 0\\
\upsilon \Sigma_r & 0 & (1-\lambda)I_r
\end{array}
\right|\\
&\!=\!(1-\omega-\lambda)^{m-r} \mathop \Pi \limits_{i=1}^r \left[(1-\lambda)(1-\omega-\lambda-\gamma \upsilon \sigma_i^2)+\omega \upsilon \sigma_i^2\right]
\end{array}
\]
by the computational property of determinant. So $\lambda$ is an eigenvalue of $\hat H_{\omega, \gamma, \upsilon}$ if and only if $\lambda =1-\omega$
or
 \[\mathop \Pi \limits_{i=1}^r [(1-\lambda)(1-\omega-\lambda-\gamma \upsilon \sigma_i^2)+\omega \upsilon \sigma_i^2]=0 \]
holds. In other words, $\lambda = 1- \omega$ or $\lambda$ is a root of the one of the equations \eqref{polynomial}. This completes the proof.
\end{pf}
\begin{lem}\label{2}
\cite{RN215Linear} Both roots of the real quadratic equation $\lambda^2-b \lambda+c=0$ are less than unity in modulus if and only if $|c|<1$ and
\(|b|<1+c\).
\end{lem}

\begin{lem}\label{3}
Under the assumption of Theorem \ref{Theorem},
$\rho(\hat H_{\omega, \gamma, \upsilon})<1$, where $\rho(\hat H_{\omega, \gamma, \upsilon})$ denotes the spectral radius of $\hat H_{\omega, \gamma, \upsilon}$ defined by \eqref{hat}, if and only if weights $\omega, \gamma, \upsilon$ satisfy \eqref{Convergence_1}.
\end{lem}
\begin{pf}
For the given weights $\omega, \gamma, \upsilon$,
by Lemma \ref{2}, all the roots of equations in \eqref{polynomial} are less than unity in modulus if and only if
\begin{equation}\label{roots}
\left\{ \begin{array}{l}
|\sigma_i^2 \upsilon (\omega-\gamma)+1-\omega |<1,\\
|\gamma \upsilon \sigma_i^2 -(2-\omega)|<2-\omega+\sigma_i^2\upsilon(\omega-\gamma),
\end{array}
\right. i=1, 2, \cdots, r
\end{equation}
holds, so by Lemma \ref{1}, $\rho (\hat H_{\omega, \gamma, \upsilon})<1$ if and only if $|1-\omega|<1$ and \eqref{roots} holds. Or equivalently,
\begin{equation}\label{aaa}
\left\{ \begin{array}{l}
0< \omega <2,\\
\sigma_i^2 \upsilon \omega>0,\\
\omega-2<\sigma_i^2 \upsilon (\omega -\gamma)<\omega,\\
\sigma_i^2 \upsilon (\omega- 2\gamma)>2(\omega-2),
\end{array}
\right. \ i=1,2,\cdots, r.
\end{equation}
Since $\sigma_1^2 \ge \sigma_2^2 \ge \cdots \ge \sigma_r^2>0$, \eqref{Convergence_1} is equivalent to \eqref{aaa}. This completes the proof.
\end{pf}

We now turn to prove Theorem \ref{Theorem}.

\noindent \textbf{Proof of Theorem \ref{Theorem}}
\begin{pf}
 Since proving that the curve sequence $\{C^k(t)\}_{k \ge 0}$, generated by the MLSPIA method, converges to a least square fitting curve of the data set $\{Q_i\}_{i=1}^m$ is equivalent to proving that the iterative sequence $ \{P^k\}_{k \ge 0}$, generated by \eqref{DIEDAI2}, converges to a solution of \(B^TBX=B^TQ\), where $B$ and $Q$ are defined by \eqref{B}, we consider the convergence of the sequence $\{P^k\}_{k \ge 0}$.

For any weights $\omega, \gamma, \upsilon$ satisfying \eqref{Convergence_1}, let $\alpha =Q-B \beta$, where $\beta$ is any solution of $B^T B X= B^T Q$, then it can be easily verified that $(\alpha^T \ \beta^T)^T$ is a solution of the equation
\begin{equation}\label{equ}
(I_{m+n}-H_{\omega, \gamma, \upsilon})(x^T \ y^T)^T=C_{\omega, \gamma, \upsilon},
\end{equation}
where $H_{\omega, \gamma, \upsilon}$ and $C_{\omega, \gamma, \upsilon}$ are defined in \eqref{H}. That is to say, the equation \eqref{equ} is consistent.

Let $\displaystyle \binom{\Lambda}{P}$ be any solution of the equation \eqref{equ}, i.e.,
\begin{equation}\label{juzhen_fangcheng}
\begin{array}{l}
(I_{m+n}-H_{\omega, \gamma, \upsilon})
\left(
\begin{array}{*{20}{c}}
\Lambda\\
P
\end{array}
\right)
= C_{\omega, \gamma, \upsilon}.
\end{array}
\end{equation}
Then,
\begin{equation}\label{Hk}
\left( \!\begin{array}{l}
\Lambda^{k}\!-\! \Lambda\\ P^{k}\!-\!P
\end{array}\right)
\!=\!H_{\omega, \gamma, \upsilon} \left( \begin{array}{l}
\Lambda^{k-1}\!- \!\Lambda\\ P^{k-1}\!-\!P
\end{array} \right)
\!=\! \cdots
\!=\!H_{\omega, \gamma, \upsilon}^{k} \left( \begin{array}{l}
\Lambda^{0}\!-\! \Lambda\\ P^{0}\!-\!P
\end{array}\right), k \ge 1.
\end{equation}
 By Lemma \ref{3}, \eqref{Convergence_1} guarantees that $\rho (\hat H_{\omega, \gamma, \upsilon})<1$ for the weights $\omega, \gamma, \upsilon$ chosen and $\hat H_{\omega, \gamma, \upsilon}$ defined by \eqref{hat}, so $\mathop {\lim} \limits_{k \to \infty }\hat H_{\omega, \gamma, \upsilon}^k=0$. By \eqref{xiangsi}, it holds that
\[
  \mathop {\lim} \limits_{k \to \infty }(W^T H_{\omega, \gamma, \upsilon}W)^k
 =\mathop {\lim} \limits_{k \to \infty }
 \left( \begin{array}{*{20}{c}}
 \hat H_{\omega, \gamma, \upsilon}^k & 0\\
 0 & I_{n-r}
 \end{array}
 \right)
 =\left( \begin{array}{*{20}{c}}
0&0\\
0&I_{n-r}
 \end{array}
 \right),
\]
where $W$ is defined by \eqref{N}.
Therefore,
\begin{equation}\label{lim_Hk}
\mathop {\lim} \limits_{k \to \infty}{H_{\omega, \gamma, \upsilon}^k}
= W \left( \begin{array}{*{20}{c}}
0 & 0\\
0 & I_{n-r}
\end{array}
\right)W^T.
\end{equation}
By \eqref{Hk}, $\lim \limits_{k \to \infty} \displaystyle \binom{\Lambda^k-\Lambda}{P^k-P}$
exists. Consequently, $ \lim \limits_{k \to \infty} \displaystyle \binom{\Lambda^k}{P^k} $ also exists. Let
$\displaystyle \binom{\Lambda^{\infty}}{P^{\infty}} $
be the limit of
$\displaystyle \binom{\Lambda^k}{P^k} $ as $k \rightarrow \infty$, i.e.,
$\mathop {\lim} \limits_{k \to \infty}
\displaystyle \binom{\Lambda^k}{P^k} = \binom{\Lambda^{\infty}}{P^{\infty}}
$, then by \eqref{Hk} and \eqref{lim_Hk},
\begin{equation}\label{Res}
\left(\!\!
\begin{array}{*{20}{c}}
\Lambda^{\infty}\\
P^{\infty}
\end{array}
\!\!\right)
\!\!=\!\!\left(\!\!
\begin{array}{*{20}{c}}
\Lambda\\
P
\end{array}
\!\!\right)\!+\!\left(\!\!
\begin{array}{*{20}{c}}
\Lambda^{\infty}\!-\!\Lambda\\
P^{\infty}\!-\!P
\end{array}
\!\!\right)
\!=\!\left(\!\!
\begin{array}{*{20}{c}}
\Lambda\\
P
\end{array}
\!\!\right)\!+\!W \left(\!\!
\begin{array}{*{20}{c}}
0 & 0\\
0 & I_{n\!-\!r}
\end{array}\!\!\right)W^T
\left(\!\!\begin{array}{*{20}{c}}
\Lambda^0\!-\!\Lambda\\
P^0\!-\!P
\end{array}
\!\!\right)
\end{equation}
holds as $k \to \infty$ in \eqref{Hk}. Since
\begin{equation}
(\!I_{m+n}\!-\!H_{\omega, \gamma, \upsilon}\!)W
\begin{pmatrix}
0 & 0\\
0 & I_{n-r}
\end{pmatrix}=
W
\begin{pmatrix}
I_{m+r} \!-\! \hat H_{\omega, \gamma, \upsilon}& 0\\
0 & 0
\end{pmatrix}\!\!\!
\begin{pmatrix}
0 & 0\\
0 & I_{n-r}
\end{pmatrix}=0
\end{equation}
follows from \eqref{xiangsi}, we have by \eqref{juzhen_fangcheng} that
\begin{equation}\label{a1}
\begin{aligned}
(I_{m+n}-H_{\omega, \gamma, \upsilon})\!\!
\begin{pmatrix}
\Lambda^{\infty}\\
P^{\infty}
\end{pmatrix}
&\!=(\!I_{m+n}\!-\!H_{\omega, \gamma, \upsilon}\!)\!\!
\begin{pmatrix}
\Lambda\\
P
\end{pmatrix}\\
&\quad+\!\!\left[(\!I_{m+n}\!-\!H_{\omega, \gamma, \upsilon}\!)W\!\!
\begin{pmatrix}
0 & \!\!0\\
0 & \!\!I_{n-r}
\end{pmatrix}\!\right]\!\!
W^T\!
\begin{pmatrix}
\Lambda^0\!-\!\Lambda\\
P^0\!-\!P
\end{pmatrix}
\\
&= C_{\omega, \gamma, \upsilon}.
\end{aligned}\!\!\!
\end{equation}
Thus, $\displaystyle \binom{\Lambda^{\infty}}{P^{\infty}}$ is also a solution of \eqref{equ}, and  $\displaystyle \binom{\Lambda^k}{P^k}$ converges to a solution of \eqref{equ}. Since $\omega \upsilon \ne 0$, it follows from \eqref{a1} that
\begin{equation*}\label{a2}
\left\{ \begin{array}{l}
\Lambda^{\infty}+BP^{\infty} = Q,\\
B^T\Lambda^{\infty}=0.
\end{array}
\right.
\end{equation*}
 By eliminating $\Lambda^{\infty}$, we obtain
\[B^TBP^{\infty}=B^TQ,\]
which means that the sequence $\{P^k\}_{k \ge 0}$ converges to a least square fitting result. This completes the proof.
\end{pf}

\indent In the following remark, we list the expressions of the absolute error and the backward error at the $k$th step for computing  $\displaystyle \binom{\Lambda^k}{P^k}$.
\begin{rmk}
The absolute error and backward error at the $k$th step are
\[
\begin{pmatrix}
\Lambda^{k}-\Lambda^{\infty}\\
P^{k}-P^{\infty }
\end{pmatrix}
= \begin{pmatrix}
U & 0\\
0 & V_r
\end{pmatrix}
\hat H_{\omega, \gamma, \upsilon}^{k}
\begin{pmatrix}
U^T & 0\\
0 & V_r^T
\end{pmatrix}\!\!
\begin{pmatrix}
\Lambda^0-\Lambda\\
P^0-P
\end{pmatrix}\!\!,\ \ \ k\ge 0,
\]
and
\[
\!\!\!\!\!\!\begin{pmatrix}
\Lambda^{k}\!-\!\Lambda^{k-1}\\
P^{k}\!-\!P^{k-1}
\end{pmatrix}
\!=\! \begin{pmatrix}
U & 0\\
0 & V_r
\!\end{pmatrix}
\!\!\left(\hat H_{\omega, \gamma, \upsilon}^{k}-\hat H_{\omega, \gamma, \upsilon}^{k-1}\right)\!\!
\begin{pmatrix}
U^T & 0\\
0 & V_r^T
\end{pmatrix}\!\!
\begin{pmatrix}
\Lambda^0\!-\!\Lambda\\
P^0\!-\!P
\end{pmatrix}\!\!,\ \ k \!\ge \!1,\!
\]
separately, where
$\displaystyle \binom{
\Lambda^{\infty}}
{P^{\infty}}$
 is defined by \eqref{Res} and $\displaystyle \binom{
\Lambda}
{P}$ is any solution of \eqref{juzhen_fangcheng}, $\hat H_{\omega, \gamma, \upsilon}$ is defined by \eqref{hat}
and $V_r$ is the first $r$ column of $V$ in \eqref{N}.

In fact, by \eqref{hat} and \eqref{Res}, it follows that
\begin{eqnarray*}
\begin{pmatrix}
\Lambda^{k}-\Lambda^{\infty}\\
P^{k}-P^{\infty}
\end{pmatrix}
\!\!\!\!&=&\!\!\!\!\begin{pmatrix}
\Lambda^{k}-\Lambda\\
P^{k}-P
\end{pmatrix}-\begin{pmatrix}
\Lambda^{\infty}-\Lambda \\
P^{\infty}-P
\end{pmatrix}\\[0.5cm]
\!\!\!\!&=&\!\!\!\! W \begin{pmatrix}
\hat H_{\omega, \gamma, \upsilon}^{k} & 0\\
0 & 0
\end{pmatrix}W^T \begin{pmatrix}
\Lambda^{0}-\Lambda \\
P^{0}-P
\end{pmatrix}\\[0.5cm]
\!\!\!\!&=&\!\!\!\! \begin{pmatrix}
U & 0\\
0 & V_r
\end{pmatrix} \hat H_{\omega, \gamma, \upsilon}^{k} \begin{pmatrix}
U^T & 0\\
0 & V_r^T
\end{pmatrix}\!\!
\begin{pmatrix}
\Lambda^{0}-\Lambda \\
P^{0}-P
\end{pmatrix},
\end{eqnarray*}
and then it holds that
\begin{eqnarray*}
\begin{pmatrix}
\Lambda^{k}\!-\!\Lambda^{k-1}\\
P^{k}\!-\!P^{k-1}
\end{pmatrix}
\!\!\!\!\!&=&\!\!\!\! \begin{pmatrix}
\Lambda^{k}-\Lambda^{\infty}\\
P^{k}-P^{\infty}
\end{pmatrix}- \begin{pmatrix}
\Lambda^{k-1}-\Lambda^{\infty}\\
P^{k-1}-P^{\infty}
\end{pmatrix}\\
\!\!\!\!\!&=&\!\!\!\! \begin{pmatrix}
U & 0\\
0 & V_r
\end{pmatrix}\!\!
\left( \hat H_{\omega, \gamma, \upsilon}^{k}\!-\!\hat H_{\omega, \gamma, \upsilon}^{k-1}\right)\!\!
\begin{pmatrix}
U^T & 0\\
0 & V_r^T
\end{pmatrix}\!\!
\begin{pmatrix}
\Lambda^0\!-\!\Lambda\\
P^0\!-\!P
\end{pmatrix}\!\!.\!\!
\end{eqnarray*}
\end{rmk}
\section{The comparison of  the LSPIA and MLSPIA methods}\label{choice}

 In this section, we compare the convergence rates of the LSPIA method and  MLSPIA method.

 Usually, the convergence rate of any method in the format similar to \eqref{juzhen_xingshi} is measured by the spectral radius of the iteration matrix. The smaller the spectral radius is, the faster the convergence rate is.

In the following, we give a comparison result.
\begin{thm}\label{Theorem_2}
Under the same assumption of Theorem \ref{Theorem}, we have
\begin{itemize}
  \item [1)] the convergence rate of the MLSPIA method \eqref{DIEDAI2} for approximating curves at the weights \(\omega=\omega^*\), \(\gamma=\gamma^*\), \(\upsilon=\upsilon^*\) is
\begin{equation}\label{less}
\rho(\hat H_{\omega^{*}, \gamma^{*}, \upsilon^{*}})=\frac{\sigma_1-\sigma_r}{\sigma_1+\sigma_r} \le \frac{\sigma_1^2-\sigma_r^2}{\sigma_1^2+\sigma_r^2}<1.
\end{equation}
Here $\omega^*, \gamma^*, \upsilon^*$ are defined by
\begin{equation}\label{para}
\omega^*=\gamma^*=\frac{4\sigma_1\sigma_r}{(\sigma_1+\sigma_r)^2},\ \upsilon^*=\frac{1}{\sigma_1\sigma_r},
\end{equation}
$\rho(\hat H_{\omega^{*}, \gamma^{*}, \upsilon^{*}})$ is the spectral radius of the matrix $\hat H_{\omega^{*}, \gamma^{*}, \upsilon^{*}}$, and $\sigma_1$ and $\sigma_r$ are the largest and the smallest singular values of the matrix $B$, respectively, as  in  Theorem \ref{Theorem}.
  \item [2)] when $\sigma_1 \ne \sigma_r$, the fastest convergence rate of the MLSPIA method for approximating curves is faster than that of the LSPIA method defined in \cite{RN183LSPIA} for the case that the NTP basis is linearly independent.
\end{itemize}
\end{thm}
\begin{pf}
By Lemma \ref{1}, the eigenvalues of the iterative matrix $H_{\omega^*, \gamma^*, \upsilon^*}$ are $1-\omega^*$ and all the roots of the corresponding quadratic polynomials defined by \eqref{polynomial} with $\omega, \gamma, \upsilon$ replaced by $\omega^*, \gamma^*, \upsilon^*$ for each $i=1, 2, \cdots, r$.

\indent Clearly, it holds that
\begin{equation} \label{modulus1}
|1-\omega^*|=\left(\frac{\sigma_1-\sigma_r}{\sigma_1+\sigma_r}\right)^2<\frac{\sigma_1-\sigma_r}{\sigma_1+\sigma_r}<1.
\end{equation}

Since $x^2-(\sigma_1^2+\sigma_r^2)x+\sigma_1^2\sigma_r^2 \le 0$ when $\sigma_r^2 \le x \le \sigma_1^2$, for each $i=1, 2, \cdots, r$, the discriminant of the corresponding quadratic polynomial defined by \eqref{polynomial}, for each $i=1,2,\cdots, r$, satisfies
\begin{equation}\begin{array}{ll}
D_i \!\!\!\!\! & \triangleq  (\sigma_i^2 \gamma^* \upsilon^*+\omega^*)^2-4\sigma_i^2\upsilon^* \omega^* \\
&= \left[\sigma_i^2 \dfrac{4\sigma_1\sigma_r}{(\sigma_1+\sigma_r)^2}\dfrac{1}{\sigma_1\sigma_r}+\dfrac{4\sigma_1\sigma_r}{(\sigma_1+\sigma_r)^2}\right]^2-
4\sigma_i^2\dfrac{1}{\sigma_1\sigma_r}\dfrac{4\sigma_1\sigma_r}{(\sigma_1+\sigma_r)^2}\\[0.4cm]
& = \dfrac{16\left[\sigma_i^4-(\sigma_1^2+\sigma_r^2)\sigma_i^2+\sigma_1^2\sigma_r^2\right]}{(\sigma_1+\sigma_r)^4} \le 0.
\end{array}
\end{equation}

If $D_i = 0$ for some $1 \le i \le r$, i.e., $\sigma_i^2 = \sigma_1^2$ or $\sigma_i^2 = \sigma_r^2$, then the corresponding quadratic polynomial has multiple roots
\begin{equation}\label{root_1}
\lambda_{i1}=\lambda_{i2} =\left\{ \begin{array}{l}
\dfrac{\sigma_r-\sigma_1}{\sigma_r+\sigma_1}, \ \ \ if\ \ \ \sigma_i=\sigma_1,\\[0.2cm]
\dfrac{\sigma_1-\sigma_r}{\sigma_r+\sigma_1}, \ \ \ if\ \ \ \sigma_i=\sigma_r
\end{array}
\right.
\end{equation}
for these $1 \le i \le r$.
In a word,
\begin{equation}\label{modulus2}
|\lambda_{i1}|=|\lambda_{i2}|=\frac{\sigma_1-\sigma_r}{\sigma_1+\sigma_r}
\end{equation}
 holds for these $1 \le i \le r$.\\
\indent If $D_i < 0$ for some $1 \le i \le r$,
then the corresponding quadratic polynomial has two conjugate imaginary roots, denoted as $\lambda_i$ and ${\overline \lambda}_i$ without inconvenience.
Then it holds that
\begin{equation}\label{modulus3}
|{\overline \lambda}_i|=|\lambda_i|=\sqrt{{\overline \lambda}_i\lambda_i} =\sqrt{1-\omega^*}=\frac{\sigma_1-\sigma_r}{\sigma_1+\sigma_r},
\end{equation}
for these $1 \le i \le r$.

It follows from \eqref{modulus1}, \eqref{modulus2} and \eqref{modulus3} that
\[
\rho(\hat H_{\omega^*, \gamma^*, \upsilon^*})= \max \left\{\frac{\sigma_1-\sigma_r}{\sigma_1+\sigma_r},\ \left(\frac{\sigma_1-\sigma_r}{\sigma_1+\sigma_r}\right)^2\right\}=\frac{\sigma_1-\sigma_r}{\sigma_1+\sigma_r}<1.
\]
By Lemma \ref{2}, the inequality $\rho (\hat H_{\omega^*, \gamma^*, \upsilon^*})<1$ means that $\omega^*, \gamma^*, \upsilon^*$ satisfy \eqref{Convergence_1}.

The fact that
\[
\frac{\sigma_1^2-\sigma_r^2}{\sigma_1^2+\sigma_r^2}-\frac{\sigma_1-\sigma_r}{\sigma_1+\sigma_r}=
\frac{2\sigma_1\sigma_r(\sigma_1-\sigma_r)}{(\sigma_1^2+\sigma_r^2)(\sigma_1+\sigma_r)}\ge 0
\]
deduces that the inequality of \eqref{less} is true. 1) is proved.

For the case when the NTP basis is linearly independent, $B$ is of full column rank, that is to say, $r=rank(B)=n$, we have $\sigma_r=\sigma_n$ and that $\sigma_1^2, \sigma_2^2, \cdots, \sigma_n^2$ are all the positive eigenvalues of $BB^T$ by \eqref{qiyi}. Since $B^T B$ and $BB^T$  have the same nonzero eigenvalues (including multiples), $\sigma_1^2, \sigma_2^2, \cdots, \sigma_n^2$ are just all the eigenvalues of the symmetric and positive definite matrix $B^T B$. It is shown in \cite{RN183LSPIA} that $\dfrac{\sigma_1^2-\sigma_n^2}{\sigma_1^2+\sigma_n^2}$ is the fastest convergence rate of the LSPIA method,
so, by the inequality in \eqref{less}, $\rho (H_{\omega^*, \gamma^*, \upsilon^*})$, the spectral radius of $H_{\omega, \gamma, \upsilon}$ at weights $\omega=\omega^*$, $\gamma=\gamma^*$, $\upsilon=\upsilon^*$, is smaller than the fastest convergence rate of the LSPIA method \cite{RN183LSPIA} if $\sigma_1 \ne \sigma_n$.
We can now conclude that the fastest convergence rate of the MLSPIA method is faster than that of the LSPIA method \cite{RN183LSPIA} for this case when $\sigma_1 \ne \sigma_n$. 2) is proved.

The proof is completed.
\end{pf}

Based on Theorem \ref{Theorem_2}, $\omega=\omega^*$, $\gamma=\gamma^*$ and $\upsilon=\upsilon^*$ may be a suitable choice to guarantee the faster convergence of the MLSPIA method.

\begin{rmk}In this paper, we do not consider the problem of how weights influence the convergent rate of the algorithm. The answer to this problem may be long due to the difficulty of obtaining the expression of the $\rho (\hat H_{\omega, \gamma, \upsilon})$ for any weights $\omega, \gamma, \upsilon$ that guarantee the convergence. It is interesting and hard and could be the subject of a new paper.
\end{rmk}

\begin{rmk}\label{remark_2}After all, it is a least square fitting, some unwanted undulations, a not well local approximation to a segment of the original curve, might happen sometimes. In this case, local adjustments of knots and increment of control points during the iterative procedure may improve the fitting result,  as is done  for the  LSPIA method in Section $5.2$ of \cite{RN183LSPIA}. Since we are only interested in accelerate the convergence rate of iteration methods here, we omit the discussion of this part here.
\end{rmk}

\section{The MLSPIA method for surfaces}\label{surface}

In this section, we will extend the MLSPIA method to fitting surfaces.

For the given data point set $\{Q_{ij}\}_{i=1, j=1}^{m_1, m_2}$ and  any initial control point set  $\{P_{ij}\}_{i=1, j=1}^{n_1, n_2}$  in the same vector space with $m_1>n_1, m_2>n_2$,  let $\{t_i, s_j\}_{i=1, j=1}^{m_1, m_2}$, satisfying $0 \le t_1<t_2< \cdots <t_{m_1} \le 1$ and $0 \le s_1<s_2< \cdots <s_{m_2} \le 1$, be two increasing sequences and $\{\phi_i(t)\}_{i=1}^{n_1}$ and $\{\psi_i(t)\}_{i=1}^{n_2}$ be two NTP bases of the vector spaces of real functions defined on $[0,1]$. The MLSPIA method approximates the surface iteratively by the tensor product surfaces
\begin{equation}\label{sur_1}
C^{k+1}(t,s)=\sum\limits_{i=1}^{n_1} \sum\limits_{j=1}^{n_2} \phi_i(t)\psi_j(s)P_{ij}^k,\ \ \ \ 0 \le t,s \le 1.
\end{equation}
Here,
\begin{equation}\label{sur_2}
P_{ij}^{k+1}=P_{ij}^{k}+\Delta_{ij}^k,\ \ \ \ i=1,2,\cdots,n_1,\ j=1,2,\cdots,n_2,
\end{equation}
and
\begin{equation}\label{sur_3}
\left\{\!\!\!
\begin{array}{lc}
\Delta_{ij}^{0}\!=\! \upsilon\! \sum \limits_{h = 1}^{m_1}\! \sum \limits_{l = 1}^{m_2} \phi_i(t_h)\psi_j(s_l) \Lambda_{hl}^0, \mathop { with\ any\ sequence}\ \{\Lambda_{hl}^0\}_{h=1,l=1}^{m_1,m_2},&\\
\Delta_{ij}^{k}\!=\!(1-\omega) \Delta_{ij}^{k-1} + \gamma \delta_{ij}^k + (\omega  - \gamma ) \delta_{ij}^{k-1},&\!\!k \ge 1,\\
\delta_{ij}^k\!=\!\upsilon \! \sum \limits_{h = 1}^{m_1}\! \sum \limits_{l = 1}^{m_2}\!\phi_i(t_h) \psi_j(s_l)\left( Q_{hl}-\! \sum \limits_{i_1 = 1}^{n_1}\! \sum \limits_{j_1 = 1}^{n_2}\!u_{i_1}(t_h) v_{j_1}(s_l)P_{i_1,j_1}^k\right),
&\!\!k \ge 0,
\end{array}
\right.\!\!\!\!\!\!
\end{equation}
for all $i = 1, 2, \cdots, m_1,\ j=1, 2, \cdots, m_2,$ where $\omega$, $\gamma$ and $\upsilon$ are three real weights.

Similar to the Theorems \ref{Theorem} and \ref{Theorem_2}, we have
\begin{thm}\label{Theorem_surface}
Let $B_1=\left( \phi_i(t_j)\right)_{m_1 \times n_1}$ and $B_2=\left( \psi_i(t_j)\right)_{m_2 \times n_2}$ with $m_1>n_1$ and $m_2>n_2$. If $rank(B_1)=r$ and $rank(B_2)=s$, and the singular values of $B_1$ and $B_2$ are $\sigma_1 \ge \sigma_2 \ge \cdots \ge \sigma_r >0$ and  $\mu_1 \ge \mu_2 \ge \cdots \ge \mu_s >0$, respectively, then
\begin{itemize}
  \item [1)] the series of the surfaces generated by the MLSPIA method defined by \eqref{sur_1}, \eqref{sur_2} and \eqref{sur_3} with weights $\omega$, $\gamma$ and $\upsilon$ converge to the least square fitting surface of  the given data point set  $(\{Q_{ij}\}_{i=1,j=1}^{m_1,m_2}\) for any  initial control  point set $\{P_{ij}^0\}_{i=1, j=1}^{n_1, n_2}$ and any  point set $\{\Lambda_{ij}^0\}_{i=1, j=1}^{m_1, m_2}$ if
\begin{equation}
\upsilon>0,\ \ \  0<\omega<2,\ \ \ \omega-\frac{\omega}{(\sigma_1\mu_1)^2 \upsilon}<\gamma<\frac{\omega}{2}-\frac{\omega-2}{(\sigma_1 \mu_1)^2 \upsilon}.
 \end{equation}
 \item [2)]the convergence rate of the MLSPIA method at the weights $\omega=\omega^*, \gamma=\gamma^*, \upsilon=\upsilon^*$ for the tensor product surfaces is
\begin{equation*}
\rho(\hat H_{\omega^{*}, \gamma^{*}, \upsilon^{*}})=\dfrac{\sigma_1\mu_1-\sigma_r\mu_s}{\sigma_1\mu_1+\sigma_r\mu_s} \le \dfrac{(\sigma_1\mu_1)^2-(\sigma_r\mu_s)^2}{(\sigma_1\mu_1)^2+(\sigma_r\mu_s)^2}<1,
\end{equation*}
Here $\omega^{*}, \gamma^{*}, \upsilon^{*}$ are defined by
\begin{equation}\label{mian_para}
\omega^{*}=\frac{4\sigma_1\mu_1\sigma_r\mu_s}{(\sigma_1\mu_1+\sigma_r\mu_s)^2},\
\gamma^{*}=\frac{4\sigma_1\mu_1\sigma_r\mu_s}{(\sigma_1\mu_1+\sigma_r\mu_s)^2},\
\upsilon^{*}=\frac{1}{\sigma_1\mu_1\sigma_r\mu_s},
\end{equation}
and $\rho(\hat H_{\omega^{*}, \gamma^{*}, \upsilon^{*}})$ is the spectral radius of the matrix $\hat H_{\omega^{*}, \gamma^{*}, \upsilon^{*}}$.
\end{itemize}
\end{thm}

The proof to the Theorem \ref{Theorem_surface} is omitted since it is similar to the case for Theorems \ref{Theorem} and \ref{Theorem_2}.

Based on $2)$ of Theorem \ref{Theorem_surface},  the wights $\omega^*, \gamma^*, \upsilon^*$, defined by  \eqref{mian_para}, may be a suitable choice to guarantee the faster convergence of the MLSPIA method for fitting the given surface.

Similar to the case  in Section \ref{choice}, how weights influence the  convergent rate of MLSPIA method for the surface fitting does not considered.

\section{Implementation and examples}\label{expetiments}

\subsection{Examples, parameter knots and weights}
In this section, five representative examples are used to test the efficiency of the MLSPIA method. They are
\begin{itemize}
\item Example 1: $205$ ($m=205$) points measured and smoothed from an airfoil-shape data.
\item Example 2: $305$ ($m=305$) points derived from a subdivision curve generated by in-center subdivision scheme.
\item Example 3: $501$ ($m=501$) points sampled uniformly from an analytic curve whose polar coordinate equation is \( r = sin \frac{\theta}{4}\ ([0,8 \pi])\).
\item Example 4: $269$ ($m=269$) points with features measured and smoothed from a G-shape font.
\item Example 5: $81 \times 81$ ($m_1=m_2=81$) points sampled from the face model.
\end{itemize}

The given data sets of Examples 1--3 and Example 5 are the same as that used in \cite{RN183LSPIA} for cubic B-spline curves fitting and  that in \cite{RN192Local} for cubic B-spline (tensor product) surface fitting, separately. The given data set used in  Example $4$ contains $269$ points which is a set smaller than that used in \cite{RN183LSPIA}. All these  data sets  are obtained from one of the authors of \cite{RN183LSPIA, RN192Local}.

In the case of curves, for the given data point set $\{Q_i\}_{i=1}^m$ and $n\le m$, as other researchers do, we use the cubic B-spline basis defined on the knot vector with the parameters, satisfying  the Sch\"{o}enberg-Whitney
condition \cite{RN401}, created by the normalized accumulated chord parameterization method. One can look for more details in $(9.5)$ and $(9.69)$ of the book \cite{RN2222QI}.

In the case of surfaces, for the given data point set $\{Q_{ij}\}_{i=1, j=1}^{m_1, m_2}$, $n_1\le m_1$ and $n_2\le m_2$, we use two groups of the cubic B-spline bases  defined separately on two knot vectors with different parameters
\begin{equation*}
\left\{\begin{array}{ll}
{\bar u_i=\dfrac{1}{m_2}  \sum \limits_{j=1}^{m_2}} t_{ij},&\ \ \ i=1, 2, \cdots, m_1,\\[2mm]
{\bar v_i=\dfrac{1}{m_1}  \sum \limits_{j=1}^{m_1}} s_{ji},&\ \ \ i=1, 2, \cdots, m_2,
\end{array}\right.
\end{equation*}
where for each $i=1, 2, \cdots, m_1$, $\{t_{ij}\}_{j=1}^{m_2}$ is generated by the normalized accumulated chord parameterization method $(9.5)$ in \cite{RN2222QI} for the data set $\{Q_{ij}\}_{j=1}^{m_2}$, and for each $i=1, 2, \cdots, m_2$, $\{s_{ji}\}_{j=1}^{m_1}$  for the data set $\{Q_{ji}\}_{j=1}^{m_1}$.

In numerical experiments, the weights appeared in the MLSPIA method are chosen as $\omega=\omega^*$, $\gamma=\gamma^*$ and $\upsilon=\upsilon^*$, where $\omega^*, \gamma^*, \upsilon^*$ are defined by \eqref{para} for Examples 1--4 and by \eqref{mian_para} for Example 5, and the weight appeared in the LSPIA method is chosen as $\mu=\mu^\ast$, where $\mu^\ast$ is the value defined by (18) in \cite{RN183LSPIA}.

\subsection{Initial control points}

The discussion on the choice of initial control points  for the MLSPIA method is not an aim of this paper. To compare the numerical experiments of the MLSPIA and LSPIA methods, only two trivial initial control points are considered.

As what is done in tests for  methods that solve systems of linear systems, a trivial choice is  that  all  initial control points  in $\{P_i^0\}_{i=1}^n$ and all points in $\{\Lambda_i^0\}_{i=1}^m$ are zero points. In this case, we can get directly via  \eqref{definition_of_P_k+1} and \eqref{Lam}, without any real numerical computation, that  $P_1^1, P_2^1, \cdots, P_n^1$ are all zero points and $\Lambda_i^1=\omega Q_i$, $i=1, 2, \cdots, m$.  To save computational cost, by replacing $P_1^0, P_2^0, \cdots, P_n^0$ with $P_1^1, P_2^1, \cdots, P_n^1$ (all zero points) and $\Lambda_1^0, \Lambda_2^0, \cdots, \Lambda_m^0$ with $\Lambda_1^1, \Lambda_2^1, \cdots, \Lambda_m^1$, we obtain the initial control points I listed below.

To compare  with  the LSPIA method, the other trivial choice is to use the same initial control points as that used by another authors in numerical experiments. Since the initial control point set is selected as a subset of the given data set in \cite{RN183LSPIA,RN181TwoWeights}, we choose  them as initial control points, too, and let $\Lambda_1^0=\Lambda_2^0=\dots=\Lambda_m^0=0$.  Thanks for the same reason to obtain  the initial control points I, we determine the initial control points II.

These two kinds of trivial initial control points for the case of curve  are

 \begin{itemize}
   \item [I:] all the $P_1^0, P_2^0, \cdots, P_n^0$ are zero  points accompanied by $\Lambda_i^0=\omega Q_i$, $i=1, 2, \cdots, m$.
   \item [II:] accompanied by $\Lambda_i^0= \omega Q_i - \omega \sum\limits_{j=1}^n \! B_j(t_i)P_j^0$, $i=1, 2, \cdots, m$, all the $P_1^0, P_2^0,$ $\cdots,$  $P_n^0$ are chosen as $Q_{f(1,m,n)}, Q_{f(2,m,n)}, \cdots, Q_{f(n,m,n)}$, where
       \begin{equation}\label{xuan_kongzhidianzhibiao}
       \left\{
       \begin{array}{l}
       f(1,m,n)=1,\\ f(i,m,n)=\left[\frac{m(i-1)}{n-1}\right]+1,\ i=2, 3, \cdots, n-1,\\ f(n,m,n)=n
       \end{array}
       \right.
       \end{equation}
       is described in the equation $(23)$ of \cite{RN183LSPIA}.
 \end{itemize}

Similarly, two kinds of trivial initial control points for the case of surfaces may be:
  \begin{itemize}
   \item [I:] all the elements of $\{P_{ij}^0\}_{i=1,j=1}^{n_1,n_2}$ are zero points accompanied by
   \[\Lambda_{i,j}^0=\omega Q_{i,j},\ \ \ \ \ \  i=1, 2, \cdots, m_1,\ j=1,2, \cdots, m_2.\]
   \item [II:] accompanied by \[\Lambda_{ij}^0= \omega Q_{ij} \!-\! \omega \sum\limits_{h=1}^{n_1} \sum\limits_{l=1}^{n_2}\! u_h(t_i)v_l(t_j)P_{hl}^0,\ \ \ i=1, 2, \cdots, m_1,\ j=1,2, \cdots, m_2,\]
        all the elements of $\{P_{ij}^0\}_{i=1,j=1}^{n_1,n_2}$ are chosen as \[P_{ij}^0=Q_{f(i,m_1,n_1),f(j,m_2,n_2)},\ \ \ i=1, 2, \cdots, n_1,\ j=1, 2, \cdots, n_2,\] where the function $f$ is defined by \eqref{xuan_kongzhidianzhibiao}. %i.e., $\{Q_i^0\}_{i=1}^m$.
 \end{itemize}

As is listed in Table \ref{sel_int_points} of Section \ref{expetiments}, we can see that under the same stop criterion, the iteration steps starting from the second kind of the initial control points are usually fewer. In the following, all the initial control points are chosen as the second kind, except in the Table \ref{sel_int_points}, where two kinds of initial control points are used.

\subsection{Numerical results }

In all of our implementations for the comparison of the  LSPIA and MLSPIA methods, the iteration process is stopped if $|E_{k}|< 10^{-7}$, where
\begin{equation}\label{wucha}
E_k= \left\|\tilde{B}^T(\tilde{B}\tilde{P}^k-\tilde{Q})\right\|_2, \ \ k \ge 0,
\end{equation}
and  the $i$th row of $\tilde{P}^k$ and $\tilde{Q}$ are coordinates of the control point $P_i^k$ at the $k$th step and corresponding $i$th ordered point $Q_i$, separately, and $\tilde{B}$ is the  collocation matrix for the cases of curves and $\tilde{B}=B_1\otimes B_2$ for the case of tensor product surfaces. Here $B_1$ and $B_2$ are defined in Theorem \ref{Theorem_surface}.

All the examples are performed on a PC with a $3.10$ GHz $64$-bit processor and $12$ GB memory via MATLAB R2014b.

In the following, firstly, we list the iteration numbers of the above examples in Table \ref{sel_int_points} for two kinds of the trivial initial control points described in Sections \ref{choice} and \ref{surface}.

In Table \ref{sel_int_points}, ``ICP" denotes the initial control points and ``NP" the number of the initial control points, ``I" and ``II" represent the first and second kinds of the initial control points separately and ``IT" is the iteration steps, $\lfloor x \rfloor$ is the biggest integer not greater than $x$. To save spaces in the table, $\lfloor x \rfloor$ is used to denote $\lfloor x \rfloor \times \lfloor x \rfloor$ for Example $5$, too.

\begin{table}[h!]
\footnotesize
\newcommand{\tabincell}[2]{\begin{tabular}{@{}#1@{}}#2\end{tabular}}
\caption{The iteration steps for the different choices of the control points}\label{sel_int_points}
\centerline{
\begin{tabular}{|c|c|c|c|c|c|c|c|c|}
\hline
&\diagbox{\scriptsize{ICP}}{\scriptsize{IT}}{\scriptsize{NP}}
& $\lfloor \frac{m}{12} \rfloor$ & $\lfloor \frac{m}{10} \rfloor$ & $\lfloor \frac{m}{8} \rfloor$ &$\lfloor \frac{m}{6} \rfloor$ & $\lfloor \frac{m}{4} \rfloor$ & $\lfloor \frac{m}{2} \rfloor$ &$\lfloor \frac{2m}{3} \rfloor$\\ \hline
\multirow{2}{*}{\scriptsize{Example $1$}}
& $\!\!\!$I$\!\!\!$& 63 & 63 & 62 & 60 & 59 & 84 & 531 \\ \cline{2-9}
& $\!\!\!$II $\!\!\!$ & 58 & 57 & 55 & 53 & 51 & 71 & 453 \\ \hline
\multirow{2}{*}{\scriptsize{Example $2$}}
& $\!\!\!$I$\!\!\!$ & 78 & 72 & 72 & 74 & 70 & 100 & 979 \\ \cline{2-9}
& $\!\!\!$II$\!\!\!$ & 69 & 64 & 62 & 61 & 54 & 75 & 946\\ \hline
\multirow{2}{*}{\scriptsize{Example $3$}}
& $\!\!\!$I$\!\!\!$ & 56 & 55 & 55 & 55 & 55 & 85 & 696 \\ \cline{2-9}
& $\!\!\!$II$\!\!\!$ & 49 & 47 & 48 & 48 & 49 & 74 & 570 \\ \hline
\multirow{2}{*}{\scriptsize{Example $4$}}
& $\!\!\!$I$\!\!\!$ & 73 & 70 & 69 & 67 & 69 & 75 & 426 \\ \cline{2-9}
& $\!\!\!$II$\!\!\!$ & 67 & 63 & 63 & 59 & 56 & 57 & 331 \\ \hline
\multirow{2}{*}{\scriptsize{Example $5$}}
& $\!\!\!$I$\!\!\!$ & 396 & 384 & 362 & 340 & 319 & 582 & 32172 \\ \cline{2-9}
& $\!\!\!$II$\!\!\!$ & 407 & 388 & 364 & 339 & 313 & 551 & 32014 \\ \hline
\multicolumn{9}{l}{\scriptsize{``ICP": the initial control points, ``NP": the number of the initial control points,}}\\
\multicolumn{9}{l}{\scriptsize{``I" and ``II": the first and second kind of the initial control points, ``IT": iteration}}\\
 \multicolumn{9}{c}{\scriptsize{steps, $\lfloor x \rfloor$: the biggest integer not greater than $x$, $\lfloor x \rfloor$ for Example $5$: $\lfloor x \rfloor \times \lfloor x \rfloor$. }}
\end{tabular}
}
\end{table}
We can see from Table \ref{sel_int_points} that the second kind of the initial control points may save iteration steps in most cases. In the following, we always use the second kind of initial control points.

Now, we begin to draw the cubic B-spline curves or cubic B-spline (tensor product) surfaces which are selected from the convergent process to the least square fitting produced by the MLSPIA method starting from the second kind of the initial control points.
% and figures of $E_k$, defined in \eqref{wucha}, and the LSPIA method, separately. with the MLSPIA method
The point sets of the above examples are fitted with $20, 30, 50, 35$, $35\times 35$ control points, respectively. In other words, $n=20, 30, 50, 35$ for Examples $1-4$ and $n_1=n_2=35$ for Example $5$. Here,  the number of the control points is chosen  as the same as that in \cite{RN183LSPIA} in Examples $1-3$, separately. Since the shapes of Examples $4-5$ are slightly complexity, the number of control points are chosen greater slightly. It is chosen stochastically as an integer near $\frac{m}{8}$ in Example 4, and as $n_1\times n_2$ in  Example 5 with $n_1=n_2$ near $\frac{m_1}{3}$ (or $\frac{m_2}{3})$. The values of the weights used here are listed in Table \ref{table_2}.

In  each figure of the first four examples, the red line is the iteration curve and  green points are control points  obtained at corresponding step, and the small blue points are the given  points that shape a dotted limit curve that is to be fitted. Since Example 5 is complex, only blue grids are drawn.

\begin{figure}[H]
\setlength{\abovecaptionskip}{-0 cm}
\setlength{\belowcaptionskip}{-0 cm}
    \centering
   \subfigure[Step 0, initial fitting.]
    {
    \setlength{\abovecaptionskip}{-0 cm}
    \setlength{\belowcaptionskip}{-0 cm}
    {\includegraphics[width=0.45\textwidth]{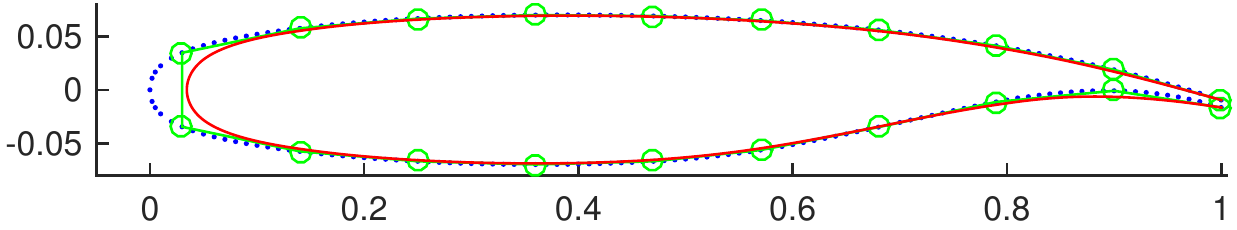}}
    \label{fig:ex1a}
    }
    \subfigure[Step 7.]
    {
    \setlength{\abovecaptionskip}{-0 cm}
    \setlength{\belowcaptionskip}{-0 cm}
    {\includegraphics[width=0.45\textwidth]{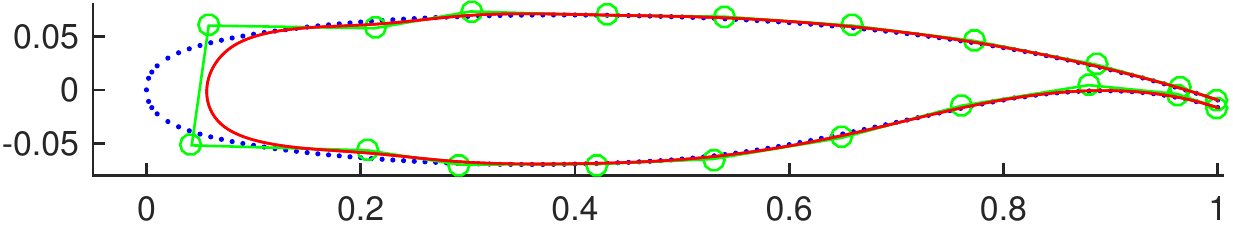}}
    \label{fig:ex1b}
    }
    \subfigure[Step 14.]
    {
    \setlength{\abovecaptionskip}{0.cm}
    \setlength{\belowcaptionskip}{-0.cm}
     {\includegraphics[width=0.45\textwidth]{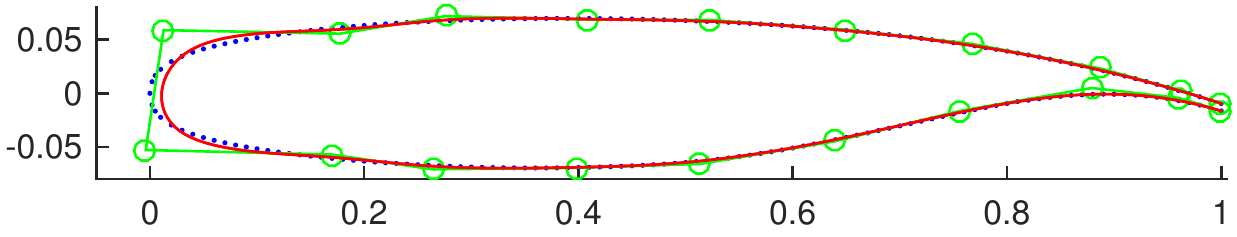}}
    \label{fig:ex1c}
    }
     \subfigure[Step 57, fitting result.]
    {
    \setlength{\abovecaptionskip}{0.cm}
    \setlength{\belowcaptionskip}{-0.cm}
    {\includegraphics[width=0.45\textwidth]{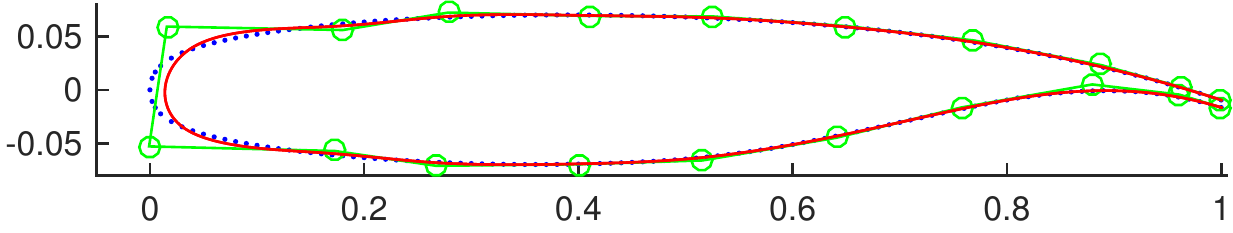}}
    \label{fig:ex1d}
    }
    \caption{An airfoil-shape data set of 205 data points is fitted by a cubic B-spline curve with 20 control points.}
    \label{fig:figure1}
\end{figure}

For Example 1, the cubic B-spline curves generated by the MLSPIA method at the initial step, $7$th and $14$th steps, and $57$th step are shown in Figure \ref{fig:figure1}.\vskip -0.3cm

\begin{figure}[H]
    \centering
    \subfigure[Initial fitting.]
    {
    {\includegraphics[width=0.2\textwidth]{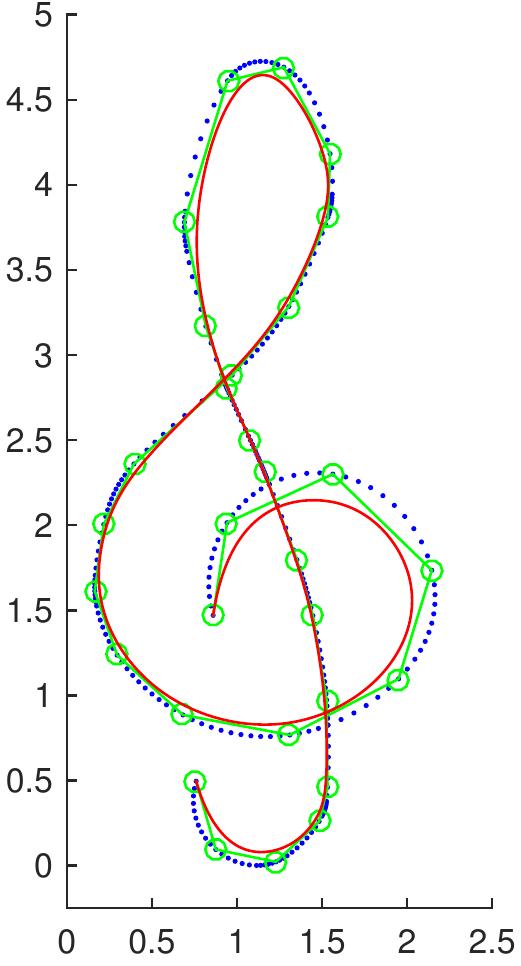}}
    \label{fig:ex2a}
    }
   \subfigure[Step 5.]
    {
    {\includegraphics[width=0.2\textwidth]{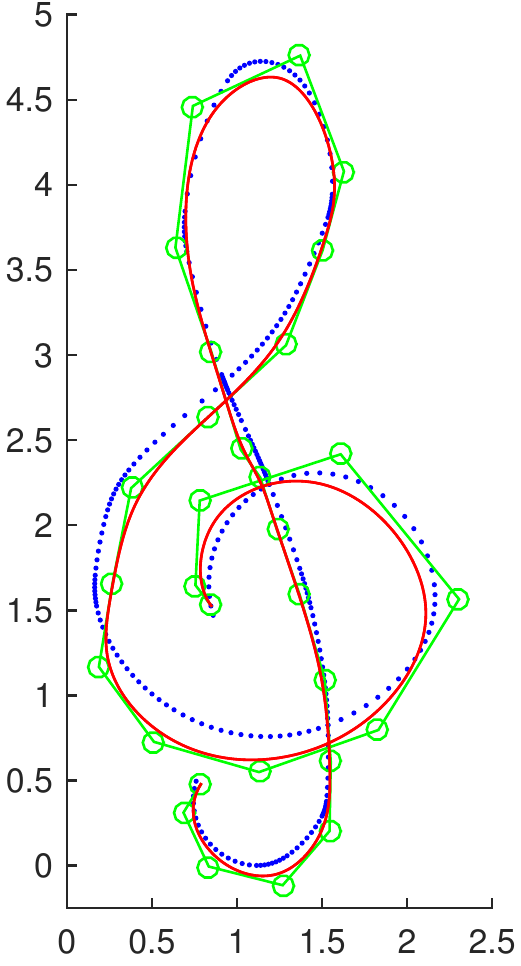}}
    \label{fig:ex2b}
    }
   \subfigure[Step 10.]
    {
    {\includegraphics[width=0.2\textwidth]{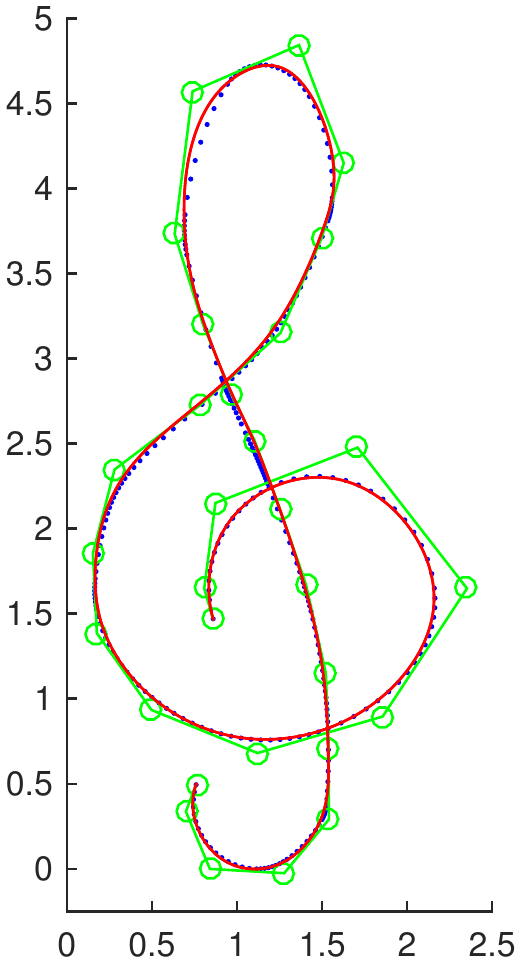}}
    \label{fig:ex2c}
    }
    \subfigure[Fitting result.]
    {
    {\includegraphics[width=0.2\textwidth]{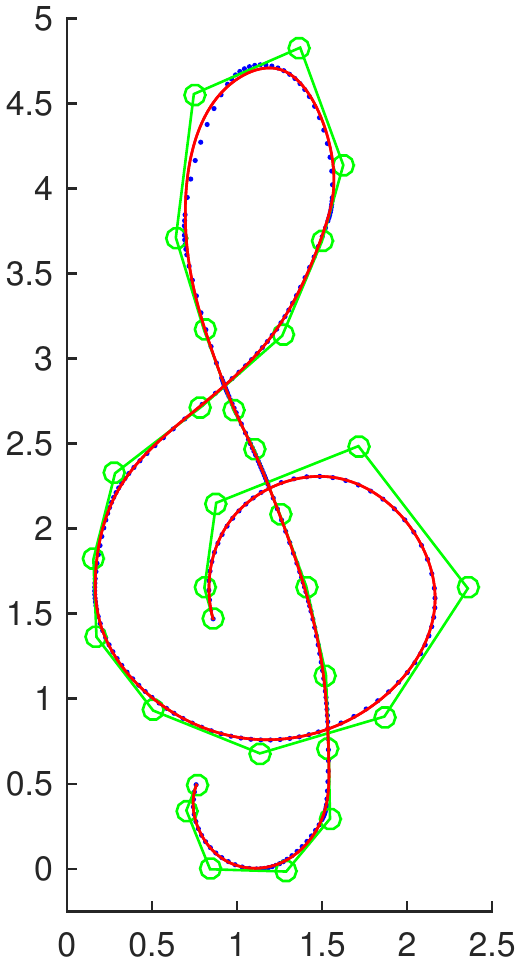}}
    \label{fig:ex2d}
    }
    \caption{A point set of 305 points is fitted by a cubic B-spline curve with 30 control points.}
    \label{fig:figure2}
 \end{figure}

For Example 2, the cubic B-spline curves generated by the MLSPIA method with the second kind of initial points at the initial step, $5$th and $10$th steps, and $64$th step are shown in Figure \ref{fig:figure2}.
 \vspace{-0.3cm}

\begin{figure}[H]
    \centering
    \subfigure[Initial fitting.]
    {
    {\includegraphics[width=0.23\textwidth]{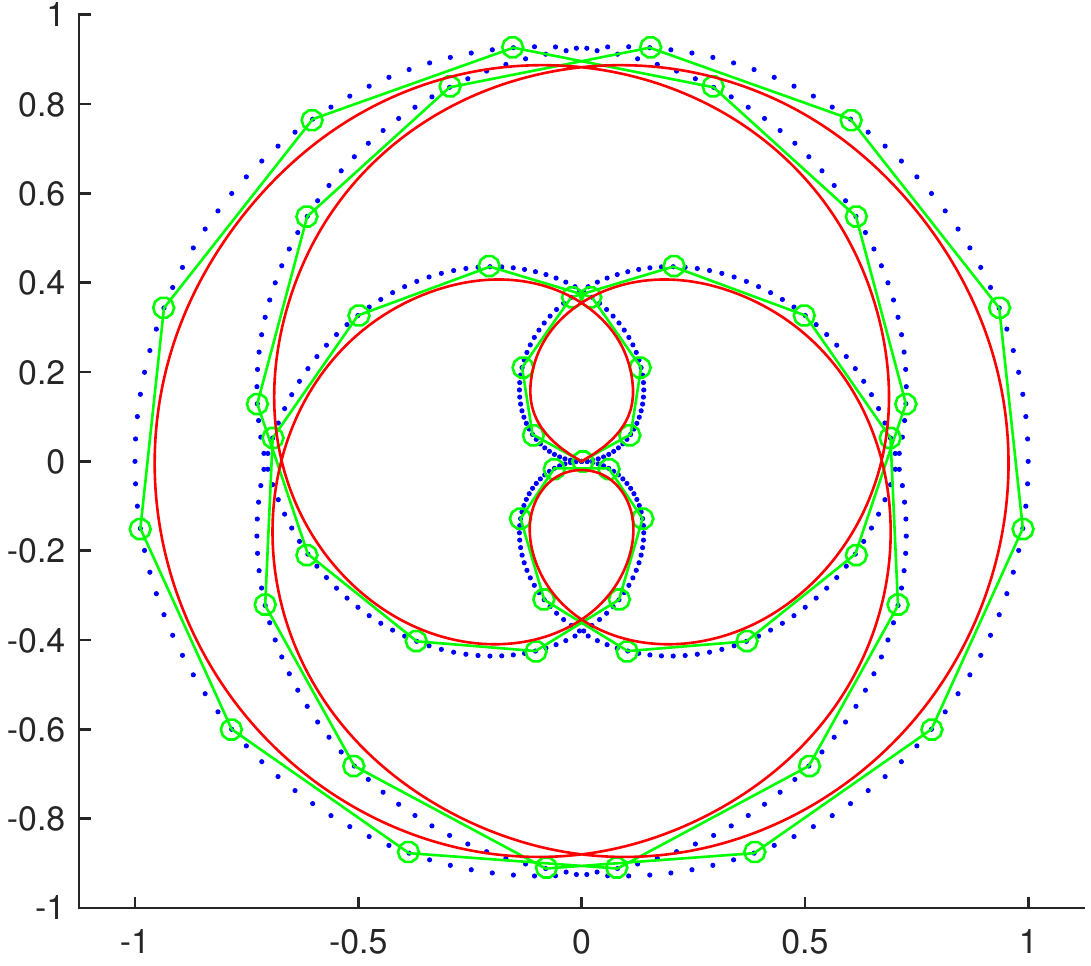}}
    \label{fig:ex3a}
    }\!\!\!\!
   \subfigure[Step 5.]
    {
     {\includegraphics[width=0.23\textwidth]{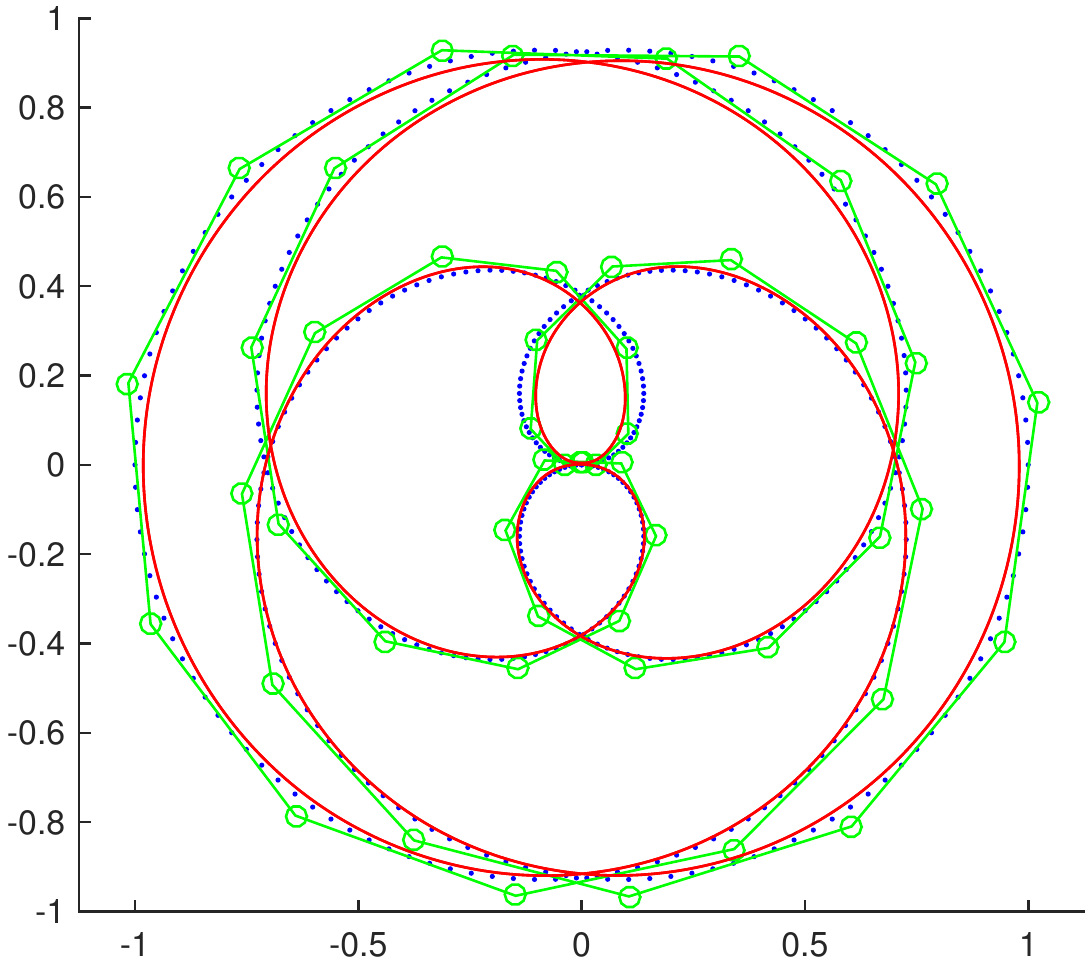}}
    \label{fig:ex3b}
    }\!\!\!\!
   \subfigure[Step 10.]
    {
    {\includegraphics[width=0.23\textwidth]{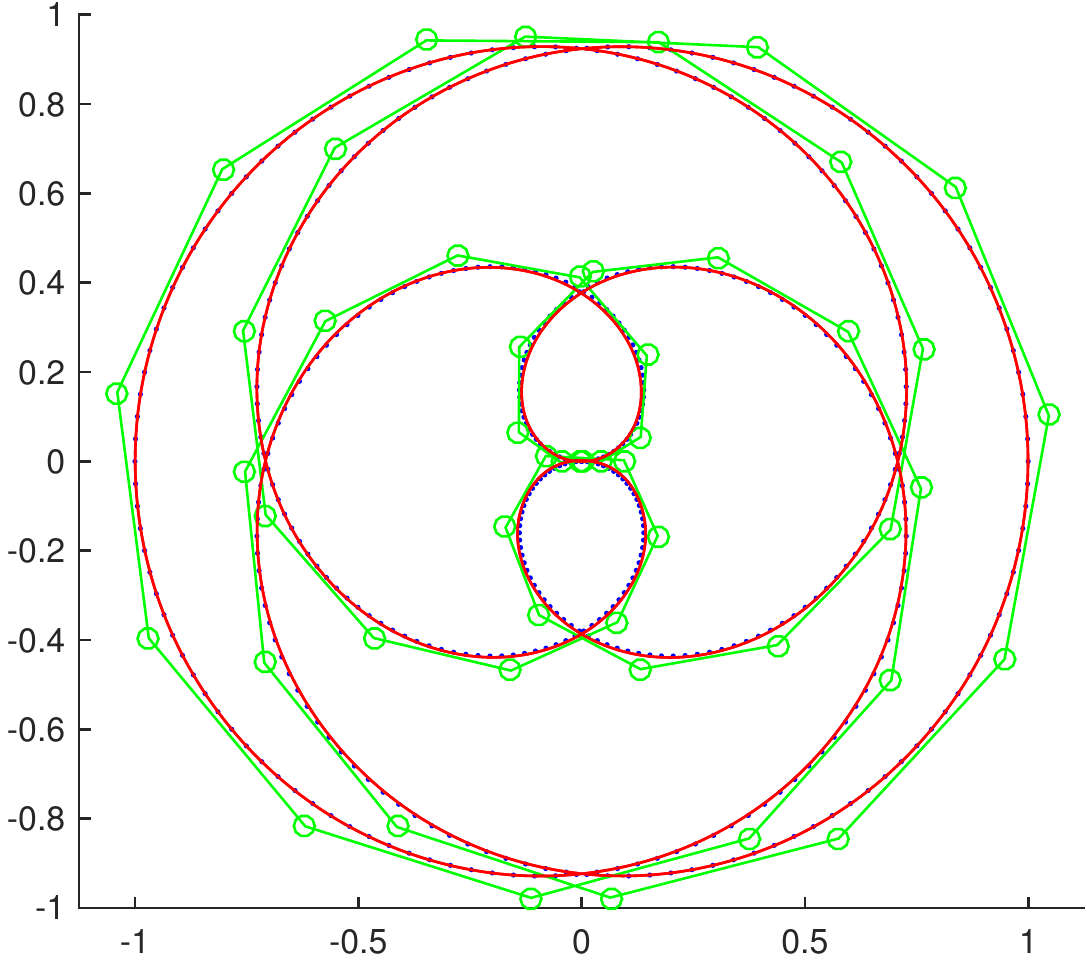}}
    \label{fig:ex3c}
    }\!\!\!\!
    \subfigure[Fitting result.]
    {
    {\includegraphics[width=0.23\textwidth]{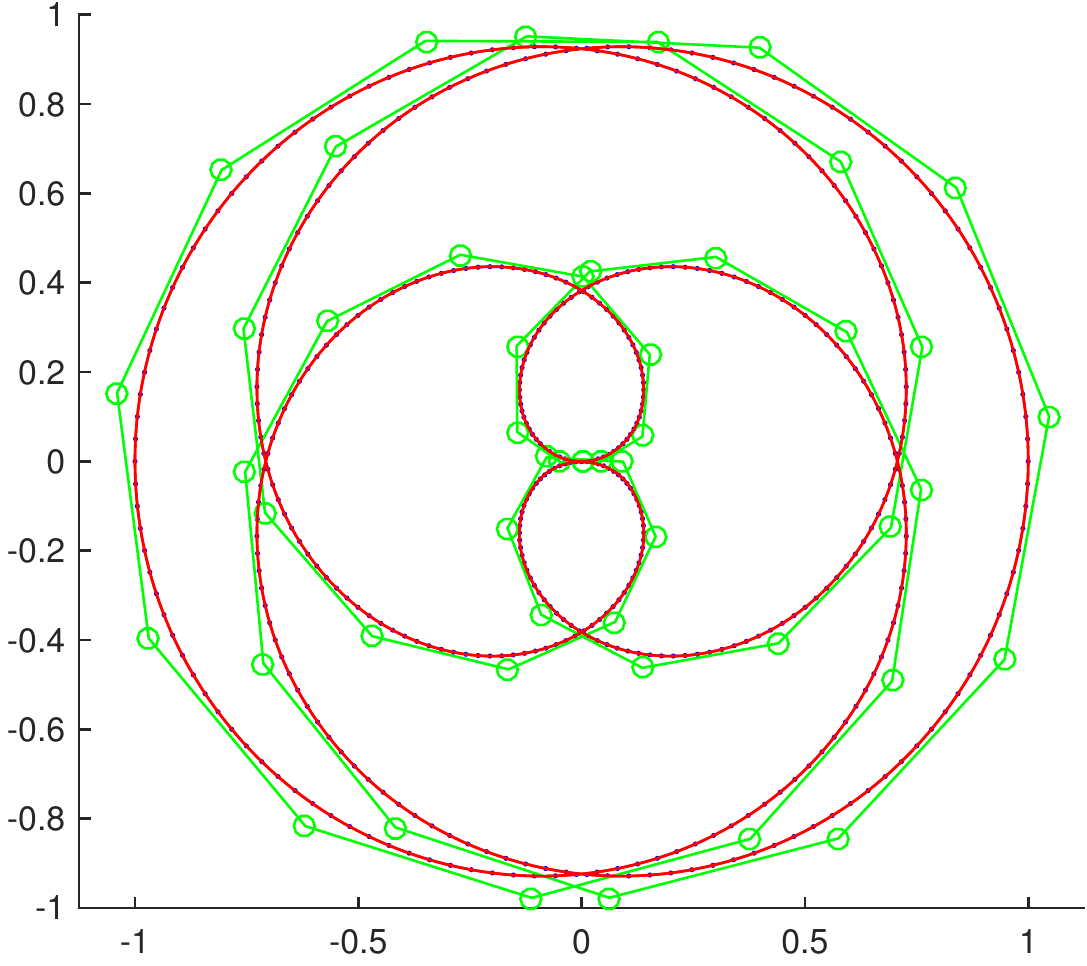}}
    \label{fig:ex3f}
    }
    \caption{A point set of 501 points is fitted by a cubic B-spline curve with 50 control points.}
    \label{fig:figure3}
\end{figure}

For Example 3, the cubic B-spline curves generated by the MLSPIA method with the second kind of initial points at the initial step, $5$th and $10$th steps, and $47$th step are shown in Figure \ref{fig:figure3}.
 %\vspace{-0.1cm}

\begin{figure}[H]
\setlength{\abovecaptionskip}{-0 cm}
\setlength{\belowcaptionskip}{-0 cm}
    \centering
   \subfigure[Initial fitting.]
    {
    \setlength{\abovecaptionskip}{-0 cm}
    \setlength{\belowcaptionskip}{-0 cm}
    {\includegraphics[width=0.23\textwidth]{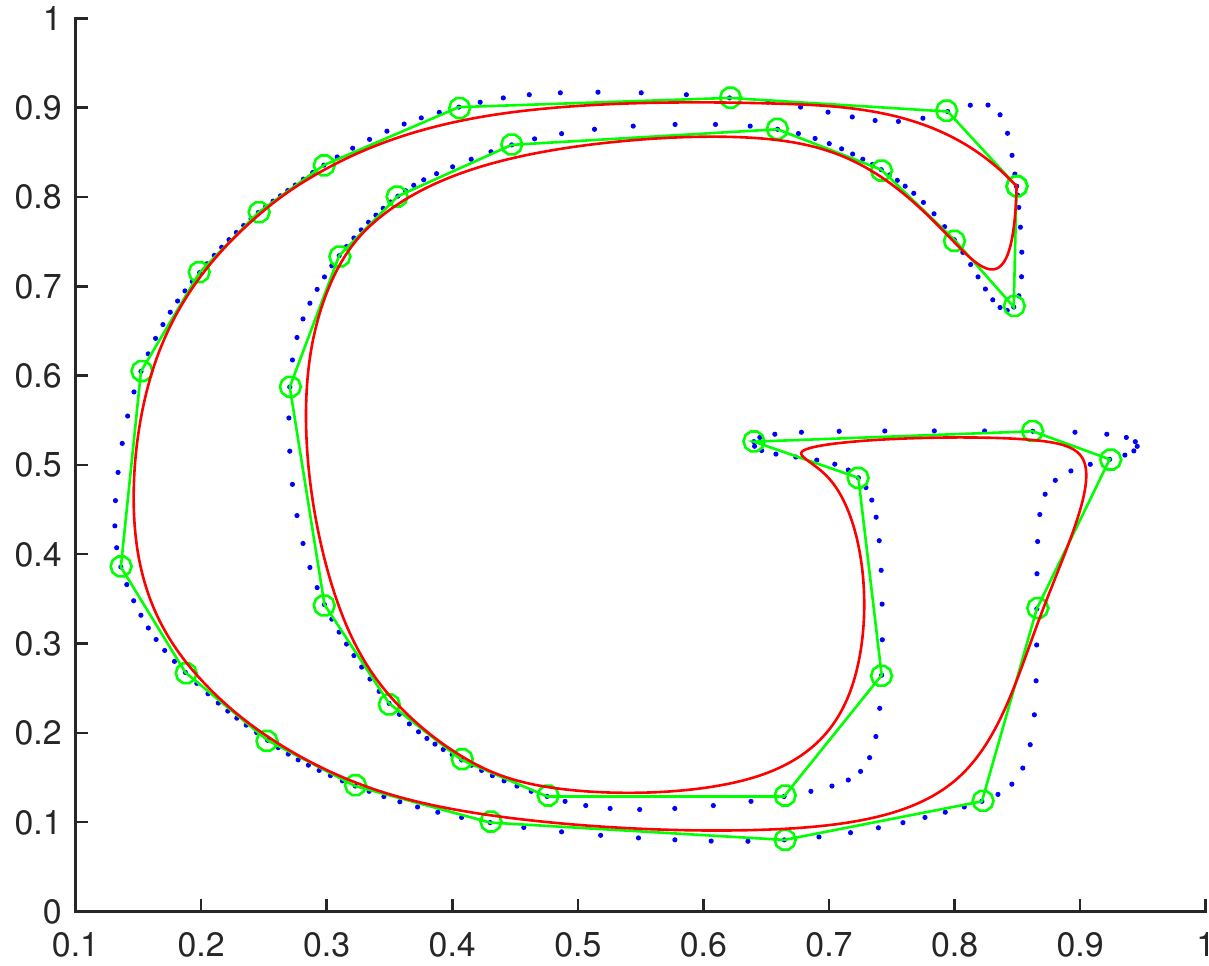}}
    \label{fig:ex1a}
    }\!\!\!\!\!\!\!\!\!
    \subfigure[Step 20.]
    {
    \setlength{\abovecaptionskip}{-0 cm}
    \setlength{\belowcaptionskip}{-0 cm}
    {\includegraphics[width=0.23\textwidth]{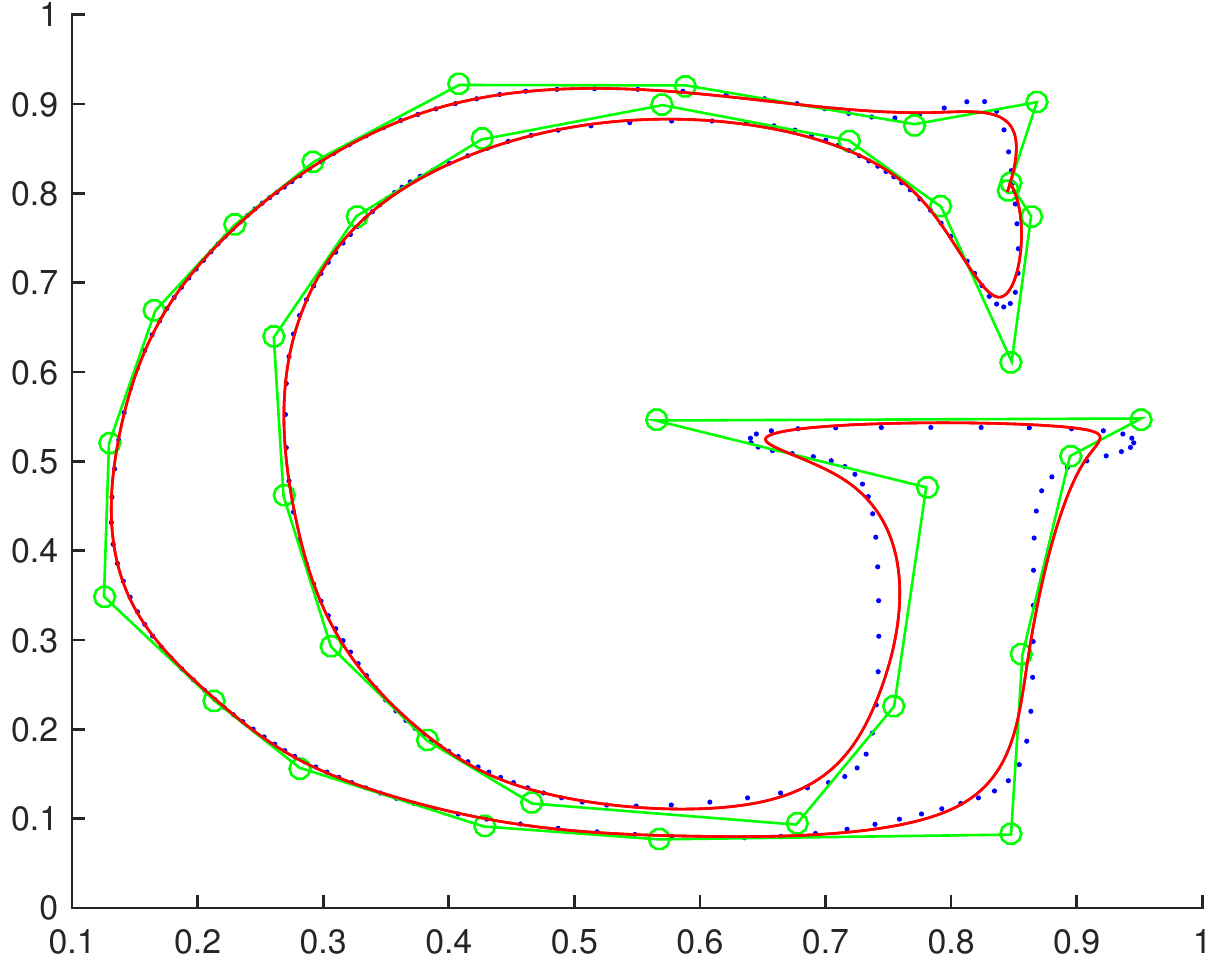}}
    \label{fig:ex1b}
    }\!\!\!\!\!\!\!\!
    \subfigure[Step 40.]
    {
    \setlength{\abovecaptionskip}{0.cm}
    \setlength{\belowcaptionskip}{-0.cm}
     {\includegraphics[width=0.23\textwidth]{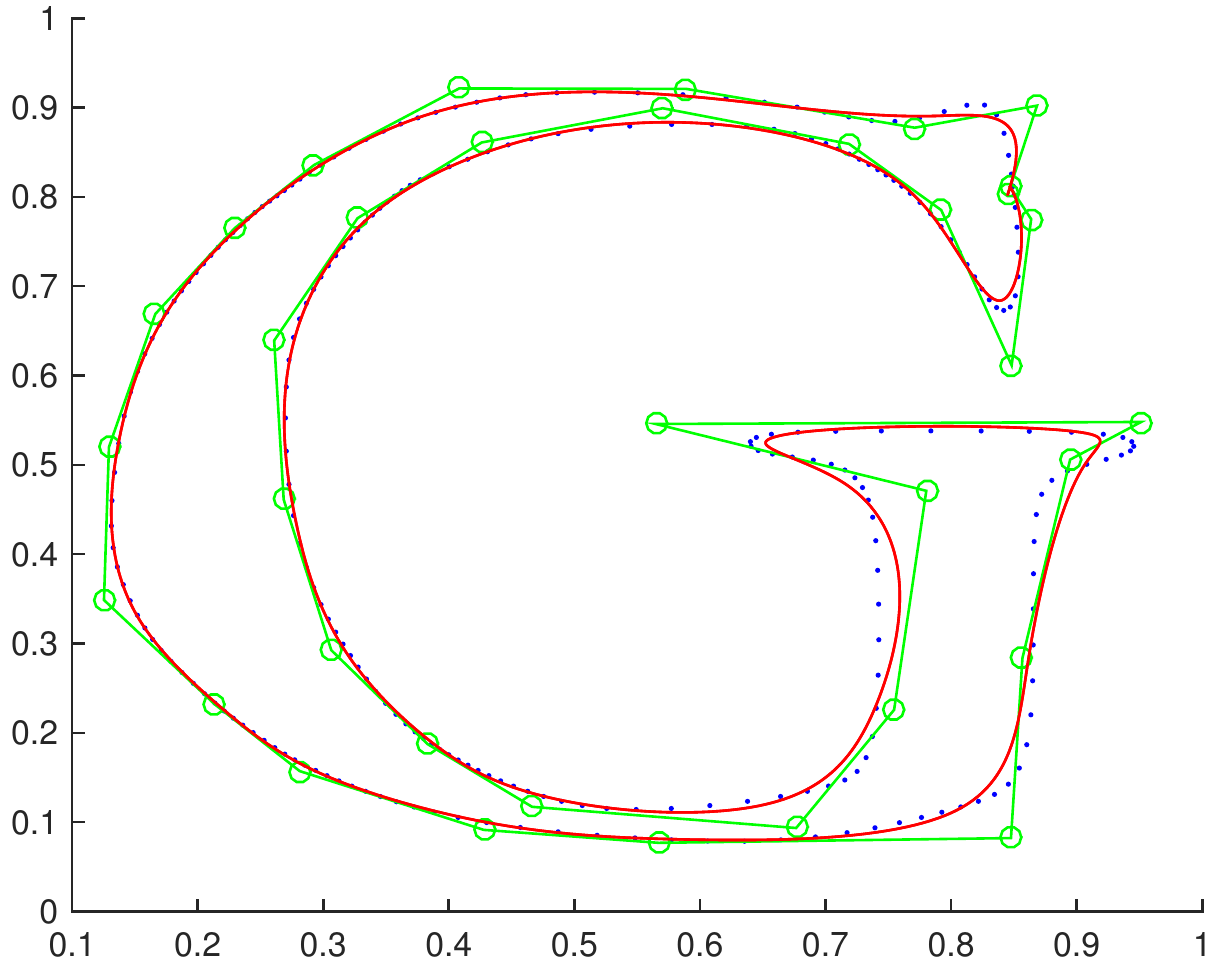}}
    \label{fig:ex1c}
    }\!\!\!\!\!\!\!\!
     \subfigure[Fitting result.]
    {
    \setlength{\abovecaptionskip}{0.cm}
    \setlength{\belowcaptionskip}{-0.cm}
    {\includegraphics[width=0.23\textwidth]{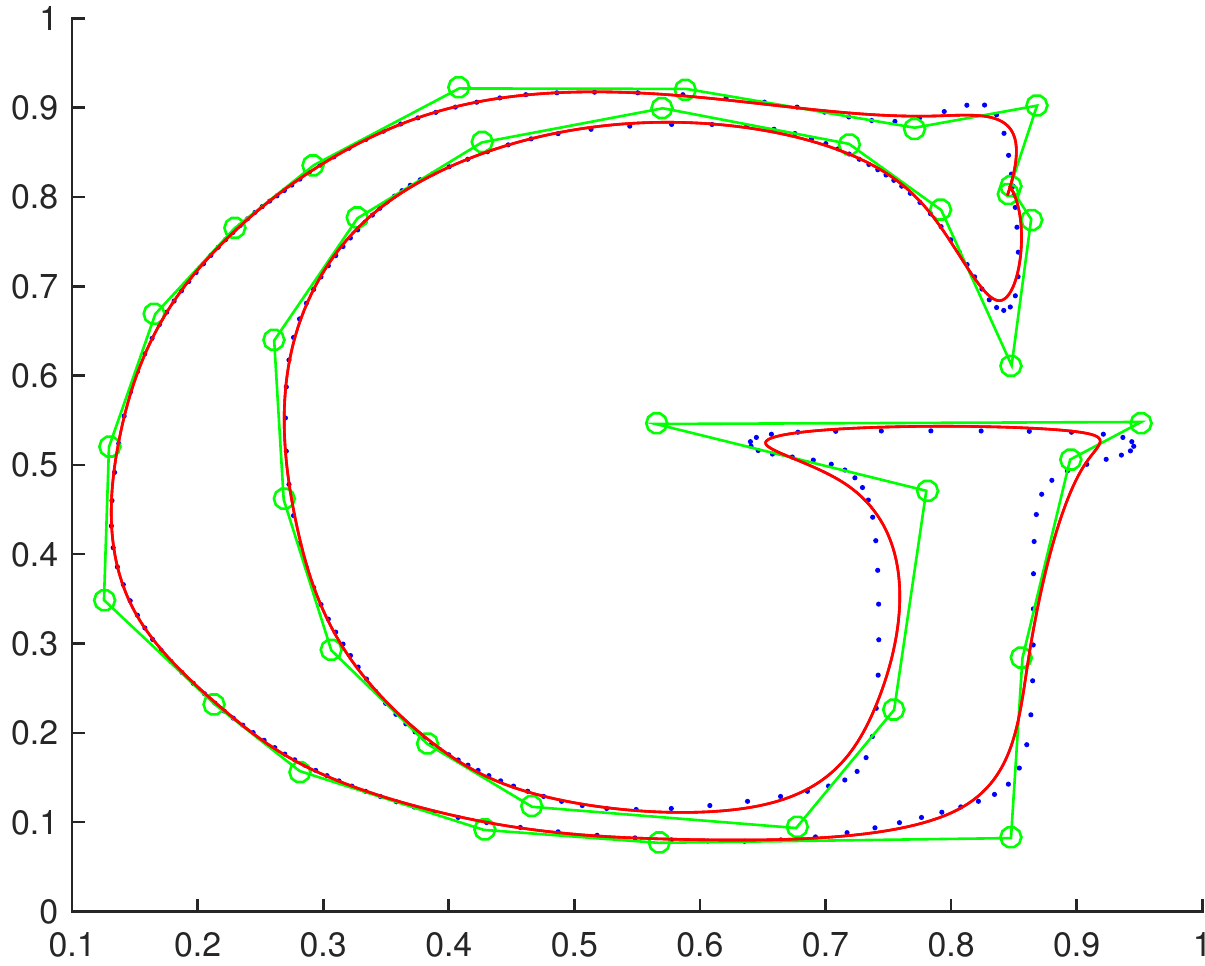}}
    \label{fig:ex1d}
    }
    \caption{A point set of 269 points is fitted by a cubic B-spline curve with 35 control points.}
    \label{fig:figure4}
\end{figure}

For Example $4$, the cubic B-spline curves generated by the MLSPIA method with the second kind of initial points at the initial step, $20$th and $40$th steps, and $63$th step are shown in Figure \ref{fig:figure4}.
%\vskip-0.5cm

\begin{figure}[H]
\setlength{\abovecaptionskip}{-0 cm}
\setlength{\belowcaptionskip}{-0 cm}
    \centering
    \subfigure[Initial surface.]
    {
    \setlength{\abovecaptionskip}{-0 cm}
    \setlength{\belowcaptionskip}{-0 cm}
    {\includegraphics[width=0.26\textwidth]{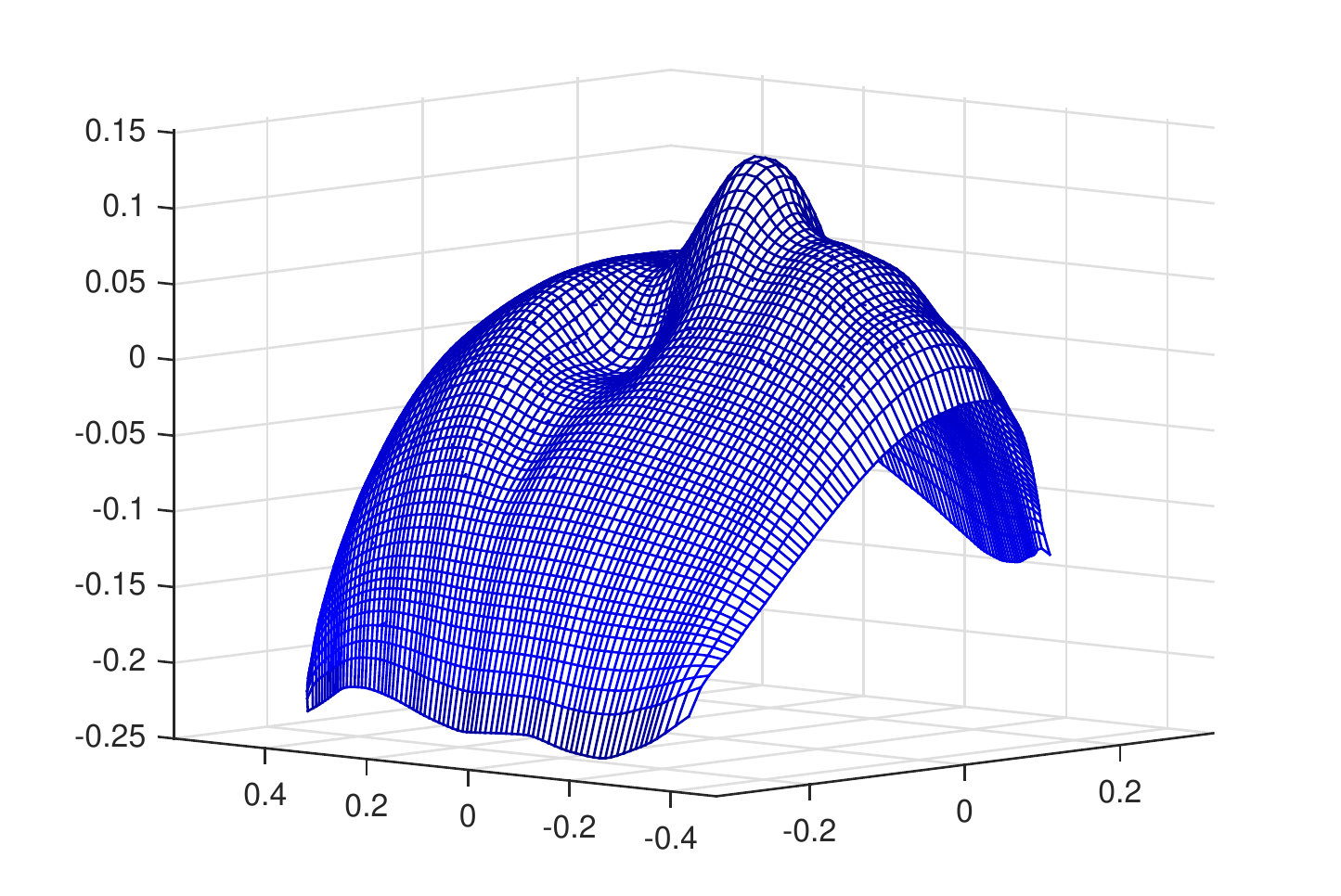}}
    \label{fig:ex1b}
    }\!\!\!\!\!\!\!\!\!\!\!\!\!\!\!
    \subfigure[Step 1000.]
    {
    \setlength{\abovecaptionskip}{0.cm}
    \setlength{\belowcaptionskip}{-0.cm}
     {\includegraphics[width=0.26\textwidth]{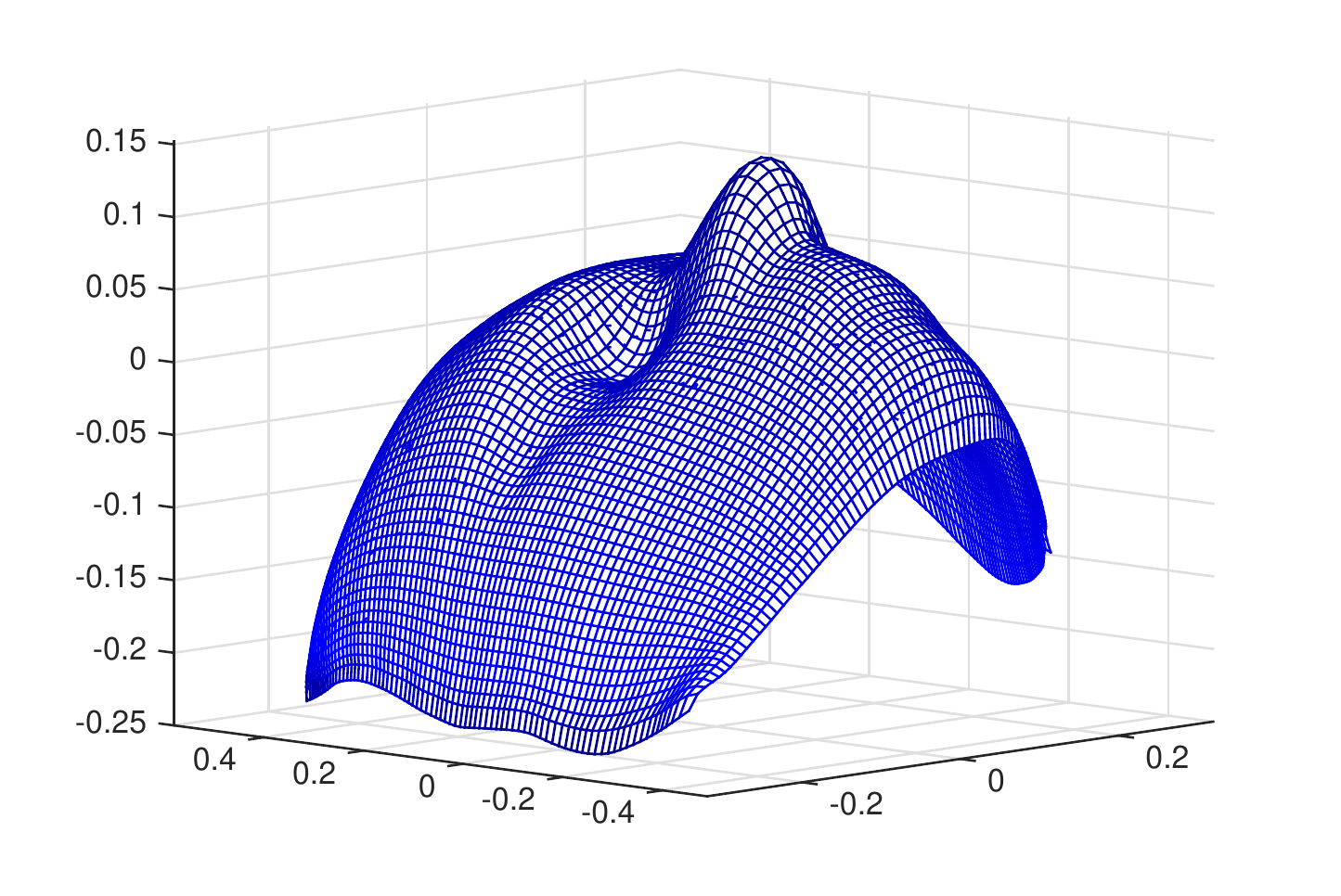}}
    \label{fig:ex1c}
    }\!\!\!\!\!\!\!\!\!\!\!\!\!\!\!\!
   \subfigure[Fitting result.]
    {
    \setlength{\abovecaptionskip}{0.cm}
    \setlength{\belowcaptionskip}{-0.cm}
     {\includegraphics[width=0.26\textwidth]{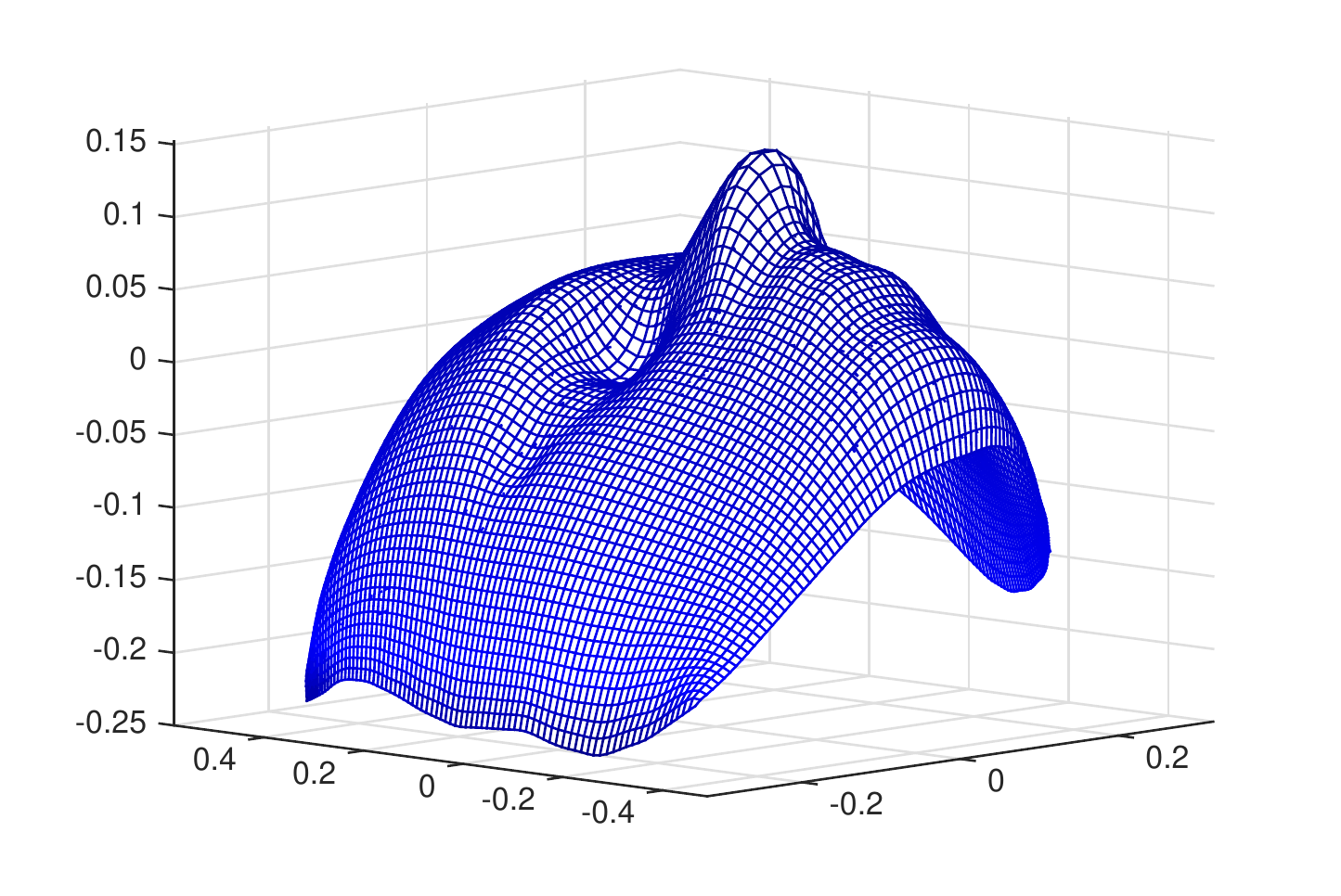}}
    \label{fig:ex1c}
    }\!\!\!\!\!\!\!\!\!\!\!\!\!\!\!\!
    \subfigure[Original surface.]
    {
    \setlength{\abovecaptionskip}{-0 cm}
    \setlength{\belowcaptionskip}{-0 cm}
    {\includegraphics[width=0.26\textwidth]{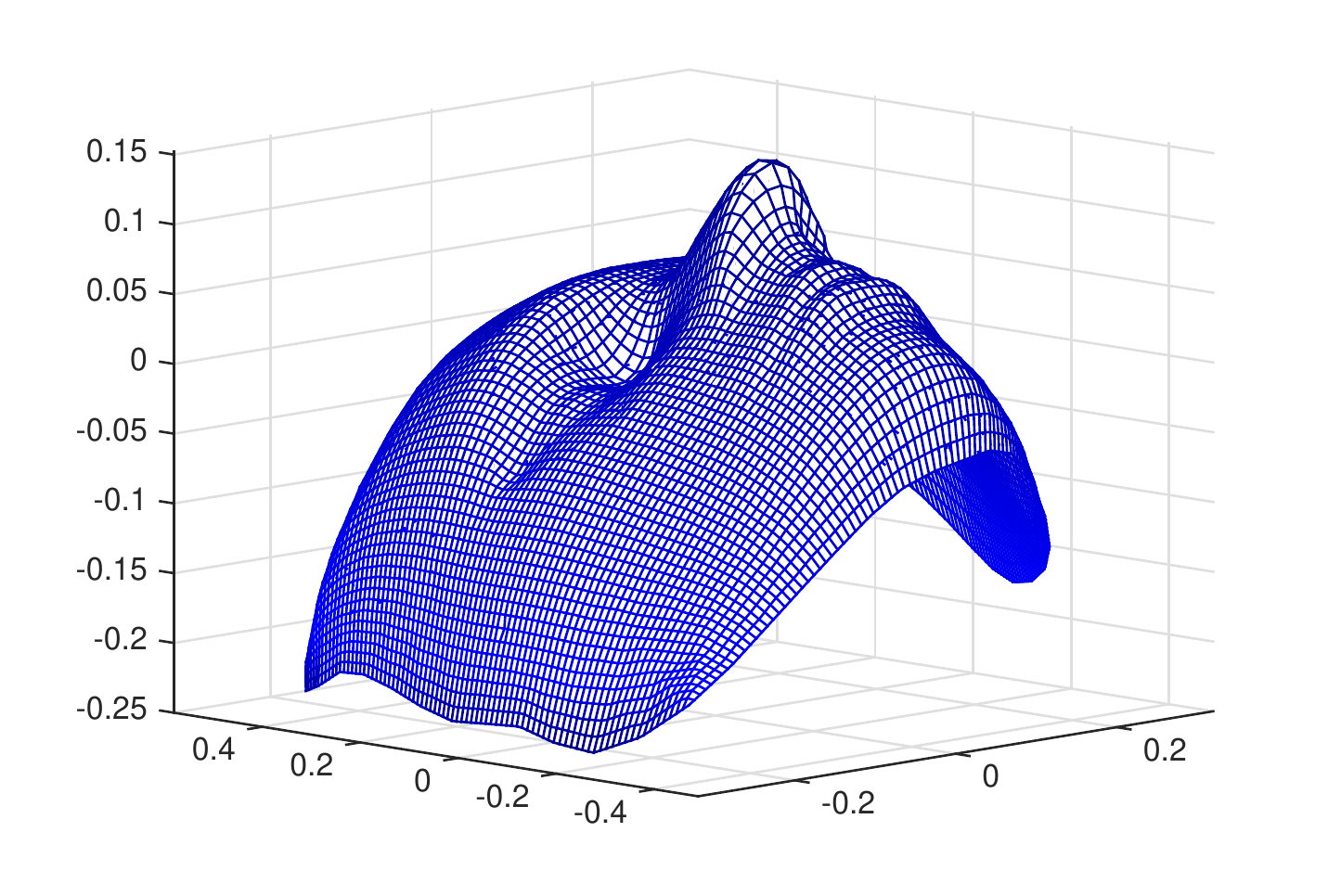}}
    \label{fig:ex1a}
    }
    \caption{A point set of $81 \times 81$ points is fitted by a cubic B-spline tensor product surface with $35 \times 35$ control points..}
    \label{fig:figure5}
\end{figure}

 For Example 5, the cubic B-spline tensor product surfaces generated by the MLSPIA method with the second kind of initial points at the initial step, $200$th and $362$th steps, and the original surface are shown in Figure \ref{fig:figure5}.

It can be seen that the fittings  at the left side in Figure \ref{fig:figure1} and around the corner at the right part in Figure \ref{fig:figure4} are not very well.  It is the same as what is stated in \cite{RN183LSPIA} for the LSPIA method. To improve the fitting result, as is described in Remark \ref{remark_2}, it is suggested that effective techniques, such as the incremental data fitting and the shape preserving ones that are stated in \cite{RN183LSPIA}, should be used in practice.

In Table \ref{table_2}, columns ``IT" (iteration steps), ``CPU time" (the average of CPU times in ten times)  and ``Weights" (the weights used in numerical experiments)  for the LSPIA and MLSPIA methods are listed. The value listed in the column ``MaxDeviations" for each example is defined by
\begin{equation*}
 \max_{0\le t\le 1}\|C_{LS}(t)-C_{MLS}(t)\|,
 \end{equation*}
 \begin{table}[h!]
\scriptsize
\centering  %
\caption{Comparison of iteration steps, CPU times and maximum deviations at chosen weights }\label{table_2}
\begin{tabular}{|c|c|c|c|c|c|}
\hline
Examples&Methods & IT & CPU time(s) & Weights & MaxDeviations \\
\hline
\multirow{3}{*}{ 1}
&\multirow{2}{*}{MLSPIA} &\multirow{2}{*}{57} &\multirow{2}{*}{0.00094} & $\omega^\ast=\gamma^\ast=0.564044574432$&\multirow{3}{*}{$5.7533\times10^{-11}$}\\
 &&&&$\nu^\ast=0.410174207181$&\\
 \cline{2-5}% \hline
& LSPIA & 224& 0.0030 &$\mu^\ast=0.161116795141$&\\ \hline
\multirow{3}{*}{2}
&\multirow{2}{*} {MLSPIA}&\multirow{2}{*}{64} & \multirow{2}{*}{0.0043} &$\omega^\ast=\gamma^\ast=0.540473573417$&\multirow{3}{*}{$8.0937\times10^{-10}$}\\
 &&& &$\nu^\ast=0.441899756336$ &\\% \hline
 \cline{2-5}
& LSPIA & 268& 0.0188 &$\mu^\ast=0.163638791357$&\\ \hline
\multirow{3}{*}{3}
&\multirow{2}{*}{MLSPIA}&\multirow{2}{*}{47} &\multirow{2}{*}{0.0045}&  $\omega^\ast=\gamma^\ast=0.589932226424$&\multirow{3}{*}{$1.1266\times10^{-9}$}\\
 % \hline
 &&&&$\nu^\ast=0.418520492365$&\\
 \cline{2-5}
 & LSPIA &156& 0.0109 &$\mu^\ast=0.175097063057$&\\ \hline
\multirow{3}{*}{4}
&\multirow{2}{*} {MLSPIA}&\multirow{2}{*}{63} & \multirow{2}{*}{0.0058} &$\omega^\ast=\gamma^\ast=0.524388806864$&\multirow{3}{*}{$2.5985\times10^{-10}$}\\ % \hline
 &&&&$\nu^\ast=0.572453504530$&\\
 \cline{2-5}
& LSPIA &272& 0.0124 &$\mu^\ast=0.203433134434$&\\ \hline
\multirow{3}{*}{5}
&\multirow{2}{*} {MLSPIA}&\multirow{2}{*}{362} & \multirow{2}{*}{5.4167} &$\omega^\ast=\gamma^\ast=0.132873535878$&\multirow{3}{*}{$1.8688\times10^{-10}$}\\ % \hline
 &&&&$\nu^\ast=4.404202175610$&\\
 \cline{2-5}
& LSPIA & 7576&130.6242&$\mu^\ast=0.313423823743$ &\\ \hline
\multicolumn{6}{l}{\scriptsize{``IT": iteration steps, ``CPU time": average CPU times in ten runs, ``Weights": $\omega^\ast, \gamma^\ast, \nu^\ast$ are}}\\
\multicolumn{6}{c}{\scriptsize{ defined by \eqref{para} for Examples 1--4 and by \eqref{mian_para} for Example 5, and $\mu^\ast$ is defined by (18) in \cite{RN183LSPIA}.}}\\
\end{tabular}
\end{table}
where $C_{LS}(t)$ and $C_{MLS}(t)$ with $0\le t\le 1$ are the least square fitting curves (or surfaces) obtained by the LSPIA and MLSPIA methods, separately. We can see that the MLSPIA method needs less iteration steps and spend less CPU time than the LSPIA method, and that the maximum deviations are very tiny.

\begin{figure}[H]
\setlength{\abovecaptionskip}{-0 cm}
\setlength{\belowcaptionskip}{-0 cm}
    \centering
    \subfigure[Examples 1-2.]
    {
    \setlength{\abovecaptionskip}{0.cm}
    \setlength{\belowcaptionskip}{-0.cm}
     {\includegraphics[width=0.450\textwidth]{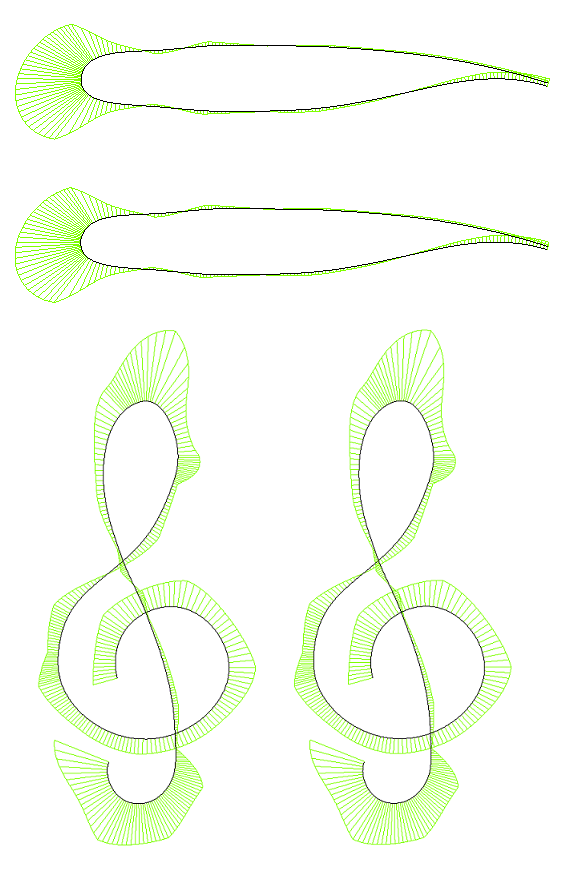}}
    \label{fig:aeroline_music_1031}
    }
    \subfigure[Examples 3-5.]
    {
    \setlength{\abovecaptionskip}{-0 cm}
    \setlength{\belowcaptionskip}{-0 cm}
    {\includegraphics[width=0.450\textwidth]{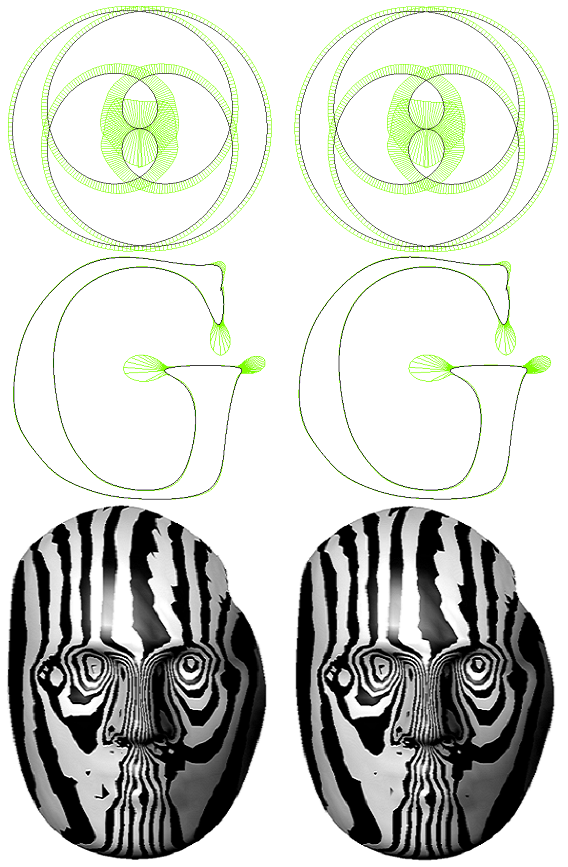}}
    \label{fig:circle_G_head}
    }
    \caption{Curvature combs for Examples 1-4 and zebra maps for Example 5: For Example 1, upper and lower figures of are the fitting results of LSPIA  and MLSPIA methods, separately. For Examples 2-5, left and right figures in each row stands for the fitting results of LSPIA  and MLSPIA methods, separately. }
    \label{fig:figure6}
\end{figure}

Curvature combs for Examples 1-4 and zebra maps for Example 5 obtained by  the LSPIA and MLSPIA methods are listed in Figure \ref{fig:figure6}, where, upper and lower figures for Example 1 are the fitting results of the LSPIA  and MLSPIA methods, separately, and left and right figures for Examples 2-5 in each row stands for the fitting results of the LSPIA  and MLSPIA methods, separately. We can see that the differences in each example are very small, as is shown in Table \ref{table_2}.

Table \ref{random_initial} lists the numerical results of the MLSPIA method in ten runs based on the sparse initial control points, obtained  by the  \textbf{rand} function of Matlab, in the rectangle $[-200000, 20000]^2$ for Examples 1-4 and in the cube $[-20000, 20000]^3$ for Example 5.  In the table, the weights are chosen as the same as that in Table \ref{table_2}, and  columns ``IT", ``CPU time" and ``Errors" list the maximum and minimum numbers of iterations, CPU times in second and errors $E_{IT}$ computed via \eqref{wucha}, separately. From the table, we can see that even for a very sparse initial control points that is far away from the given data set, the MLSPIA method will converges to the least square fitting in not too long a time.
\vskip-0.2cm
\begin{table}[h!]
\scriptsize
\centering  %
\caption{ MLSPIA method based on random initial control points}\label{random_initial}
\begin{tabular}{|c|c|c|c|c|}
\hline
Examples&IT & CPU time(s) & Errors $E_{IT}$ &Weights\\
\hline
\multirow{2}{*}{ 1}
&\multirow{2}{*}{$86\thicksim 95$} &\multirow{2}{*}{$0.0015\thicksim 0.0028$} &$6.7986\times10^{-9}\thicksim $ & $\omega^\ast=\gamma^\ast=0.564044574432$\\
&&&$9.8944\times10^{-9}$&$\nu^\ast=0.410174207181$\\
\hline
\multirow{2}{*}{2}
&\multirow{2}{*} {$92\thicksim 100$}&\multirow{2}{*}{$0.0045\thicksim 0.0072$} & $6.9669\times10^{-9}\thicksim$ &$\omega^\ast=\gamma^\ast=0.540473573417$\\
 &&&$9.8188\times10^{-9}$ &$\nu^\ast=0.441899756336$\\% \hline
\hline
\multirow{2}{*}{3}
&\multirow{2}{*}{$82\thicksim 87$}&\multirow{2}{*}{$0.0056\thicksim 0.0089$} &$7.1887\times10^{-9}\thicksim$&  $\omega^\ast=\gamma^\ast=0.589932226424$\\
 % \hline
 &&&$9,6780\times10^{-9}$&$\nu^\ast=0.418520492365$\\
\hline
\multirow{2}{*}{4}
&\multirow{2}{*} {$96\thicksim 104$}&\multirow{2}{*}{$0.0054\thicksim 0.0117$} & $7.4159\times10^{-9}\thicksim$ &$\omega^\ast=\gamma^\ast=0.524388806864$\\ % \hline
 &&&$9.7002\times10^{-9}$&$\nu^\ast=0.572453504530$\\
\hline
\multirow{2}{*}{5}
&\multirow{2}{*} {$533\thicksim 558$}&\multirow{2}{*}{$8.2565\thicksim 8.8496$} & $9.3723\times10^{-9}\thicksim$ &$\omega^\ast=\gamma^\ast=0.132873535878$\\ % \hline
 &&&$9.8929\times10^{-9}$&$\nu^\ast=4.404202175610$\\
\hline
\multicolumn{5}{l}{\scriptsize{``IT": iteration steps, ``CPU time":  CPU times, ``Weights": see Table \ref{table_2}, $``E_{IT}"$: computed via   \eqref{wucha}}.}\\
\end{tabular}
\end{table}

\section{Conclusions}
\indent In this paper, a progressive and iterative approximation method with memory for least square fitting (MLSPIA) is developed. The method constructs a series of fitting curves (surfaces) with three weights by adjusting the control points and the weighted sums iteratively. It is proved that these curves (surfaces) under the suitable choices of weights will converge to the least square fitting result even when the collocation matrix is singular, and that the fastest convergence rate of the method is faster than that of the LSPIA method proposed in \cite{RN183LSPIA} under the same condition. % Numerical experiments show that  the MLSPIA method converges faster than the LSPIA method.

\section*{Acknowledgements}
Authors would like to thank anonymous reviewers very much for their valuable comments that improve the quality of the paper, and thank Professor Lin Hongwei very much for his sending us data sets of all examples. This work is  supported by NSFC (No. 11471285 and No. 11871430).

%\section*{References}

\bibliography{Manuscipt_20190817_without_color}

\begin{thebibliography}{36}
\providecommand{\natexlab}[1]{#1}
\providecommand{\url}[1]{\texttt{#1}}
\providecommand{\href}[2]{#2}
\providecommand{\path}[1]{#1}
\providecommand{\eprint}[1]{\href{http://arxiv.org/abs/#1}{\path{#1}}}
\providecommand{\DOIprefix}{doi:}
\providecommand{\ArXivprefix}{arXiv:}
\providecommand{\URLprefix}{URL: }
\providecommand{\Pubmedprefix}{pmid:}
\providecommand{\doi}[1]{\href{http://dx.doi.org/#1}{\path{#1}}}
\providecommand{\Pubmed}[1]{\href{pmid:#1}{\path{#1}}}
\providecommand{\BIBand}{and}
\providecommand{\bibinfo}[2]{#2}
\ifx\xfnm\undefined \def\xfnm[#1]{\unskip,\space#1}\fi
%Type = Book
\bibitem[{{Piegl L, Tiller W}\!\!(1997)}]{RN2222QI}
\bibinfo{author}{{Piegl L, Tiller W}\!\!\xfnm[]}.
\newblock \bibinfo{title}{The {NURBS} book}.
\newblock \bibinfo{edition}{2nd} ed.; \bibinfo{address}{New York, USA}:
  \bibinfo{publisher}{Springer-Verlag}; \bibinfo{year}{1997}.
%Type = Book
\bibitem[{{Farin G, Hoschek J, Kim M-S}\!\!(2002)}]{Farin}
\bibinfo{author}{{Farin G, Hoschek J, Kim M-S}\!\!\xfnm[]}.
\newblock \bibinfo{title}{Handbook of computer aided geometric design}.
\newblock \bibinfo{edition}{1st} ed.; \bibinfo{publisher}{North-Holland};
  \bibinfo{year}{2002}.
%Type = Article
\bibitem[{{Borges CF, Pastva T}\!\!(2002)}]{Borgescf}
\bibinfo{author}{{Borges CF, Pastva T}\!\!\xfnm[]}.
\newblock \bibinfo{title}{Total least squares fitting of bezier and b-spline
  curves to ordered data}.
\newblock \bibinfo{journal}{Comput-Aided Des}
  \bibinfo{year}{2002};\bibinfo{volume}{19}:\bibinfo{pages}{275--289}.
%Type = Article
\bibitem[{{Wang W, Pottmann H, Liu Y}\!\!(2006)}]{RN402}
\bibinfo{author}{{Wang W, Pottmann H, Liu Y}\!\!\xfnm[]}.
\newblock \bibinfo{title}{Fitting b-spline curves to point clouds by
  curvature-based squared distance minimization}.
\newblock \bibinfo{journal}{ACM Transactions on Graphics}
  \bibinfo{year}{2006};\bibinfo{volume}{25}:\bibinfo{pages}{214--238}.
%Type = Article
\bibitem[{{Kineri Y, Wang M, Lin H, Maekawa T}\!\!(2012)}]{RN214GI}
\bibinfo{author}{{Kineri Y, Wang M, Lin H, Maekawa T}\!\!\xfnm[]}.
\newblock \bibinfo{title}{{B}-spline surface fitting by iterative geometric
  interpolation/approximation algorithms}.
\newblock \bibinfo{journal}{Comput-Aided Des}
  \bibinfo{year}{2012};\bibinfo{volume}{44}(\bibinfo{number}{7}):\bibinfo{pages}{697--708}.
%Type = Article
\bibitem[{{Lin H-W, Zhang Z-Y}\!\!(2013)}]{RN209Tspline}
\bibinfo{author}{{Lin H-W, Zhang Z-Y}\!\!\xfnm[]}.
\newblock \bibinfo{title}{An efficient method for fitting large data sets using
  {T}-splines}.
\newblock \bibinfo{journal}{SIAM J Sci Comput}
  \bibinfo{year}{2013};\bibinfo{volume}{35}(\bibinfo{number}{6}):\bibinfo{pages}{A3052--A3068}.
%Type = Article
\bibitem[{{Deng C-Y, Lin H-W}\!\!(2014)}]{RN183LSPIA}
\bibinfo{author}{{Deng C-Y, Lin H-W}\!\!\xfnm[]}.
\newblock \bibinfo{title}{Progressive and iterative approximation for least
  squares {B}-spline curve and surface fitting}.
\newblock \bibinfo{journal}{Comput-Aided Des}
  \bibinfo{year}{2014};\bibinfo{volume}{47}:\bibinfo{pages}{32--44}.
%Type = Article
\bibitem[{{Galveza A, Iglesiasa A, Avilaa A, Oteroc C, Ariasc R, Manchadoca
  C}\!\!(2015)}]{GalvezaA}
\bibinfo{author}{{Galveza A, Iglesiasa A, Avilaa A, Oteroc C, Ariasc R,
  Manchadoca C}\!\!\xfnm[]}.
\newblock \bibinfo{title}{Elitist clonal selection algorithm for optimal choice
  of free knots in b-spline data fitting}.
\newblock \bibinfo{journal}{Applied Soft Computing}
  \bibinfo{year}{2015};\bibinfo{volume}{26}:\bibinfo{pages}{90--106}.
%Type = Article
\bibitem[{{Zhang L, Ge X-Y, Tan J-Q}\!\!(2016)}]{RN181TwoWeights}
\bibinfo{author}{{Zhang L, Ge X-Y, Tan J-Q}\!\!\xfnm[]}.
\newblock \bibinfo{title}{Least square geometric iterative fitting method for
  generalized {B}-spline curves with two different kinds of weights}.
\newblock \bibinfo{journal}{Vis Comput}
  \bibinfo{year}{2016};\bibinfo{volume}{32}(\bibinfo{number}{9}):\bibinfo{pages}{1109--1120}.
%Type = Article
\bibitem[{{Lin H, Cao Q, Zhang X}\!\!(2017)}]{RN209singular}
\bibinfo{author}{{Lin H, Cao Q, Zhang X}\!\!\xfnm[]}.
\newblock \bibinfo{title}{The convergence of least-squares progressive
  iterative approximation with singular iterative matrix.}
\newblock \bibinfo{journal}{arXiv} \bibinfo{year}{2017};\URLprefix
  \url{arXiv:1707.09109.}
%Type = Article
\bibitem[{{Ebrahimi A, Loghmani BG}\!\!(2018)}]{Ebrahimi}
\bibinfo{author}{{Ebrahimi A, Loghmani BG}\!\!\xfnm[]}.
\newblock \bibinfo{title}{Shape modeling based on specifying the initial
  b-spline curve and scaled bfgs optimization method}.
\newblock \bibinfo{journal}{Multimed Tools Appl}
  \bibinfo{year}{2018};\bibinfo{volume}{77}:\bibinfo{pages}{30331--30351}.
%Type = Article
\bibitem[{{Lin H, Maekawa T, Deng C}\!\!(2018)}]{RN230Survey}
\bibinfo{author}{{Lin H, Maekawa T, Deng C}\!\!\xfnm[]}.
\newblock \bibinfo{title}{Survey on geometric iterative methods and their
  applications.}
\newblock \bibinfo{journal}{Comput-Aided Des}
  \bibinfo{year}{2018};\bibinfo{volume}{95}:\bibinfo{pages}{40--51}.
%Type = Article
\bibitem[{{Liu M-Z, Li B-J, Guo Q-J, Zhu C-G, Hu P, Shao
  Y-H}\!\!(2018)}]{RN220LIU}
\bibinfo{author}{{Liu M-Z, Li B-J, Guo Q-J, Zhu C-G, Hu P, Shao
  Y-H}\!\!\xfnm[]}.
\newblock \bibinfo{title}{Progressive iterative approximation for regularized
  least square bivariate {B}-spline surface fitting}.
\newblock \bibinfo{journal}{J Comput Appl Math}
  \bibinfo{year}{2018};\bibinfo{volume}{327}:\bibinfo{pages}{175--187}.
%Type = Article
\bibitem[{{Vaitkus M, Varady T}\!\!(2018)}]{Vaitkus}
\bibinfo{author}{{Vaitkus M, Varady T}\!\!\xfnm[]}.
\newblock \bibinfo{title}{Parameterizing and extending trimmed regions for
  tensor-product surface fitting}.
\newblock \bibinfo{journal}{Comput-Aided Des}
  \bibinfo{year}{2018};\bibinfo{volume}{104}:\bibinfo{pages}{125--140}.
%Type = Article
\bibitem[{{Lin H-W, Bao H-J, Wang G-J}\!\!(2005)}]{RN218Totally}
\bibinfo{author}{{Lin H-W, Bao H-J, Wang G-J}\!\!\xfnm[]}.
\newblock \bibinfo{title}{Totally positive bases and progressive iteration
  approximation}.
\newblock \bibinfo{journal}{Comput Math Appl}
  \bibinfo{year}{2005};\bibinfo{volume}{50}(\bibinfo{number}{3-4}):\bibinfo{pages}{575--586}.
%Type = Article
\bibitem[{{Lin H-W}\!\!(2010{\natexlab{a}})}]{RN191}
\bibinfo{author}{{Lin H-W}\!\!\xfnm[]}.
\newblock \bibinfo{title}{The convergence of the geometric interpolation
  algorithm}.
\newblock \bibinfo{journal}{Comput-Aided Des}
  \bibinfo{year}{2010}{\natexlab{a}};\bibinfo{volume}{42}(\bibinfo{number}{6}):\bibinfo{pages}{505--508}.
%Type = Article
\bibitem[{{Maekawa T, Matsumoto Y, Namiki K}\!\!(2007)}]{RN212GI}
\bibinfo{author}{{Maekawa T, Matsumoto Y, Namiki K}\!\!\xfnm[]}.
\newblock \bibinfo{title}{Interpolation by geometric algorithm}.
\newblock \bibinfo{journal}{Comput-Aided Des}
  \bibinfo{year}{2007};\bibinfo{volume}{39}(\bibinfo{number}{4}):\bibinfo{pages}{313--323}.
%Type = Article
\bibitem[{{De Boor C}\!\!(1979)}]{RN209de}
\bibinfo{author}{{De Boor C}\!\!\xfnm[]}.
\newblock \bibinfo{title}{How does {A}gee's smoothing method work?}
\newblock \bibinfo{journal}{Proceedings of the 1979 army numerical analysis and
  computers conference, ARO report}
  \bibinfo{year}{1979};\bibinfo{volume}{79-3}:\bibinfo{pages}{299--302}.
%Type = Article
\bibitem[{{Qi D, Tian Z, Zhang Y, Feng J}\!\!(1975)}]{RN208QI}
\bibinfo{author}{{Qi D, Tian Z, Zhang Y, Feng J}\!\!\xfnm[]}.
\newblock \bibinfo{title}{The method of numeric polish in curve fitting}.
\newblock \bibinfo{journal}{Acta Mathematica Sinica}
  \bibinfo{year}{1975};\bibinfo{volume}{18}:\bibinfo{pages}{173--184}.
%Type = Article
\bibitem[{{Yamaguchi F}\!\!(1977)}]{RN1977Y}
\bibinfo{author}{{Yamaguchi F}\!\!\xfnm[]}.
\newblock \bibinfo{title}{A method of designing free surfaces by computer
  display (1st report).}
\newblock \bibinfo{journal}{Precision Machinery}
  \bibinfo{year}{1977};\bibinfo{volume}{43}(\bibinfo{number}{2}):\bibinfo{pages}{168--173}.
%Type = Article
\bibitem[{{Carnicer JM, Garc\'ia-Esnaola M, Pe\~na
  JM}\!\!(1996)}]{RN225Convexity}
\bibinfo{author}{{Carnicer JM, Garc\'ia-Esnaola M, Pe\~na JM}\!\!\xfnm[]}.
\newblock \bibinfo{title}{Convexity of rational curves and total positivity}.
\newblock \bibinfo{journal}{J Comput Appl Math}
  \bibinfo{year}{1996};\bibinfo{volume}{71}(\bibinfo{number}{2}):\bibinfo{pages}{365--382}.
%Type = Article
\bibitem[{{Delgado J, Pe\~na, JM}\!\!(2003)}]{RN226Delgado}
\bibinfo{author}{{Delgado J, Pe\~na, JM}\!\!\xfnm[]}.
\newblock \bibinfo{title}{A shape preserving representation with a evaluation
  algorithm of linear complexity}.
\newblock \bibinfo{journal}{Comput Aided Geom Design}
  \bibinfo{year}{2003};\bibinfo{volume}{20}(\bibinfo{number}{1}):\bibinfo{pages}{1--10}.
%Type = Article
\bibitem[{{Lin H-W, Wang G-J, Dong C-S}\!\!(2004)}]{RN199Constructing}
\bibinfo{author}{{Lin H-W, Wang G-J, Dong C-S}\!\!\xfnm[]}.
\newblock \bibinfo{title}{Constructing iterative non-uniform {B}-spline curve
  and surface to fit data points}.
\newblock \bibinfo{journal}{Sci China Ser F}
  \bibinfo{year}{2004};\bibinfo{volume}{47}(\bibinfo{number}{3}):\bibinfo{pages}{315--331}.
%Type = Article
\bibitem[{{Lu L-Z}\!\!(2010)}]{RN193Weighted}
\bibinfo{author}{{Lu L-Z}\!\!\xfnm[]}.
\newblock \bibinfo{title}{Weighted progressive iteration approximation and
  convergence analysis}.
\newblock \bibinfo{journal}{Comput Aided Geom Design}
  \bibinfo{year}{2010};\bibinfo{volume}{27}(\bibinfo{number}{2}):\bibinfo{pages}{129--137}.
%Type = Article
\bibitem[{{Lin H-W, Zhang Z-Y}\!\!(2011)}]{RN188extended}
\bibinfo{author}{{Lin H-W, Zhang Z-Y}\!\!\xfnm[]}.
\newblock \bibinfo{title}{An extended iterative format for the
  progressive-iteration approximation}.
\newblock \bibinfo{journal}{Comput Graph}
  \bibinfo{year}{2011};\bibinfo{volume}{35}(\bibinfo{number}{5}):\bibinfo{pages}{967--975}.
%Type = Article
\bibitem[{{Shi L-M, Wang R-H}\!\!(2006)}]{RN211nurbs}
\bibinfo{author}{{Shi L-M, Wang R-H}\!\!\xfnm[]}.
\newblock \bibinfo{title}{An iterative algorithm of {NURBS} interpolation and
  approximation}.
\newblock \bibinfo{journal}{J Math Res Exposition}
  \bibinfo{year}{2006};\bibinfo{volume}{26}(\bibinfo{number}{4}):\bibinfo{pages}{735--743}.
%Type = Article
\bibitem[{{Lin H-W}\!\!(2010{\natexlab{b}})}]{RN192Local}
\bibinfo{author}{{Lin H-W}\!\!\xfnm[]}.
\newblock \bibinfo{title}{Local progressive-iterative approximation format for
  blending curves and patches}.
\newblock \bibinfo{journal}{Comput Aided Geom Design}
  \bibinfo{year}{2010}{\natexlab{b}};\bibinfo{volume}{27}(\bibinfo{number}{4}):\bibinfo{pages}{322--339}.
%Type = Article
\bibitem[{{Chen J, Wang G-J}\!\!(2011)}]{RN189triangular}
\bibinfo{author}{{Chen J, Wang G-J}\!\!\xfnm[]}.
\newblock \bibinfo{title}{Progressive iterative approximation for triangular
  {B}\'ezier surfaces}.
\newblock \bibinfo{journal}{Comput-Aided Des}
  \bibinfo{year}{2011};\bibinfo{volume}{43}(\bibinfo{number}{8}):\bibinfo{pages}{889--895}.
%Type = Article
\bibitem[{{Hu Q-Q}\!\!(2013)}]{RN221HU}
\bibinfo{author}{{Hu Q-Q}\!\!\xfnm[]}.
\newblock \bibinfo{title}{An iterative algorithm for polynomial approximation
  of rational triangular {B}\'ezier surfaces}.
\newblock \bibinfo{journal}{Appl Math Comput}
  \bibinfo{year}{2013};\bibinfo{volume}{219}(\bibinfo{number}{17}):\bibinfo{pages}{9308--9316}.
%Type = Article
\bibitem[{{Cheng F-H, Fan F-T, Lai S-H, Huang C-L, Wang J-X, Yong
  J-H}\!\!(2009)}]{RN202Loop}
\bibinfo{author}{{Cheng F-H, Fan F-T, Lai S-H, Huang C-L, Wang J-X, Yong
  J-H}\!\!\xfnm[]}.
\newblock \bibinfo{title}{Loop subdivision surface based progressive
  interpolation}.
\newblock \bibinfo{journal}{J Comput Sci Tech}
  \bibinfo{year}{2009};\bibinfo{volume}{24}(\bibinfo{number}{1}):\bibinfo{pages}{39--46}.
%Type = Article
\bibitem[{{Deng C-Y, Ma W-Y}\!\!(2012)}]{RN201WeightedPI}
\bibinfo{author}{{Deng C-Y, Ma W-Y}\!\!\xfnm[]}.
\newblock \bibinfo{title}{Weighted progressive interpolation of {L}oop
  subdivision surfaces}.
\newblock \bibinfo{journal}{Comput-Aided Des}
  \bibinfo{year}{2012};\bibinfo{volume}{44}(\bibinfo{number}{5}):\bibinfo{pages}{424--431}.
%Type = Article
\bibitem[{{Chen Z-X, Luo X-N, Tan L, Ye B-H, Chen
  J-P}\!\!(2008)}]{RN207Catmull-Clark}
\bibinfo{author}{{Chen Z-X, Luo X-N, Tan L, Ye B-H, Chen J-P}\!\!\xfnm[]}.
\newblock \bibinfo{title}{Progressive interpolation based on {C}atmull-{C}lark
  subdivision surfaces}.
\newblock \bibinfo{journal}{Computer Graphics Forum}
  \bibinfo{year}{2008};\bibinfo{volume}{27}(\bibinfo{number}{7}):\bibinfo{pages}{1823--1827}.
%Type = Article
\bibitem[{{Fan F-T, Cheng F-H, Lai S-H}\!\!(2008)}]{RN210shape}
\bibinfo{author}{{Fan F-T, Cheng F-H, Lai S-H}\!\!\xfnm[]}.
\newblock \bibinfo{title}{Subdivision based interpolation with shape control}.
\newblock \bibinfo{journal}{Computer-Aided Design and Applications}
  \bibinfo{year}{2008};\bibinfo{volume}{5}(\bibinfo{number}{1-4}):\bibinfo{pages}{539--547}.
%Type = Book
\bibitem[{{De Boor C, Rice JR}\!\!(1968)}]{RN400}
\bibinfo{author}{{De Boor C, Rice JR}\!\!\xfnm[]}.
\newblock \bibinfo{title}{Least squares cubic spline approximation {I}-Fixed
  knots}.
\newblock \bibinfo{publisher}{Computer Sciences, Purdue University};
  \bibinfo{year}{1968}.
%Type = Book
\bibitem[{{Young DM}\!\!(1971)}]{RN215Linear}
\bibinfo{author}{{Young DM}\!\!\xfnm[]}.
\newblock \bibinfo{title}{Iterative solution of large linear systems}.
\newblock \bibinfo{address}{New York-London}: \bibinfo{publisher}{Academic
  press}; \bibinfo{year}{1971}.
%Type = Article
\bibitem[{{Sch\"{o}enberg IJ, Whitney A}\!\!(1953)}]{RN401}
\bibinfo{author}{{Sch\"{o}enberg IJ, Whitney A}\!\!\xfnm[]}.
\newblock \bibinfo{title}{On polya frequency functions}.
\newblock \bibinfo{journal}{Trans Amer Math Soc}
  \bibinfo{year}{1953};\bibinfo{volume}{74}:\bibinfo{pages}{246--259}.

\end{thebibliography}

\end{document}